\documentclass[11pt,a4paper]{article}

\usepackage{booktabs}
\usepackage{amsfonts}
\usepackage{amsmath,amsthm}
\usepackage{graphicx}
\usepackage{algorithmic}
\usepackage[ruled]{algorithm2e}
\SetAlFnt{\small}
\SetAlCapFnt{\small}
\SetAlCapNameFnt{\small}
\SetAlCapHSkip{0pt}
\IncMargin{-\parindent}
\usepackage{tikz,pgfplots}
\usetikzlibrary{math}
\usepackage{xspace}
\usepackage[hidelinks]{hyperref}
\usepackage[left=2.5cm,right=2.5cm,top=3cm,bottom=3cm]{geometry}
\usepackage{cleveref}
\Crefname{ALC@unique}{Line}{Lines}

\newcounter{myalg}
\AtBeginEnvironment{algorithmic}{\refstepcounter{myalg}}
\makeatletter
\@addtoreset{ALC@unique}{myalg}
\makeatother

\newenvironment{acknowledgements}{\paragraph{Acknowledgements}}{}
\usepackage{amsopn}

\newtheorem{theorem}{Theorem}
\newtheorem{proposition}{Proposition}

\theoremstyle{remark}
\newtheorem{remark}{Remark}[section]

\newcommand{\Hsign}{S}
\newcommand{\Hhat}{\hat{H}}
\newcommand{\mybar}[1]{\makebox[0pt]{$\phantom{#1}\overline{\phantom{#1}}$}#1}
\newcommand{\tmpvec}{x}

\newcommand{\yambo}{Yambo\xspace}
\newcommand{\shao}{\textsf{shao}\xspace}
\newcommand{\gruning}{\textsf{gruning}\xspace}
\newcommand{\projectedbse}{\textsf{projectedbse}\xspace}
\newcommand{\nhe}{\textsf{non-Hermitian}\xspace}

\newcommand{\revisedtext}[2]{#2}
\newenvironment{revisedblock}{}{}
\theoremstyle{definition}
\newtheorem{definition}{Definition}[section]

\title{Variants of thick-restart Lanczos for the Bethe--Salpeter eigenvalue problem\thanks{This work was supported by grants PID2022-139568NB-I00 and RED2022-134176-T funded by MCIN/AEI/10.13039/501100011033 and by ERDF/EU. Innovation Study ISOLV-BSE has received funding through the Inno4scale project, which is funded by the European High-Performance Computing Joint Undertaking (JU) under Grant Agreement No 101118139. The JU receives support from the European Union's Horizon Europe Programme. The second author was also supported by Universitat Politècnica de València in its PAID-01-23 program.}
}

\author{Fernando Alvarruiz\thanks{D.\ Sistemes Inform\`atics i Computaci\'o, Universitat Polit\`ecnica de Val\`encia, Val\`encia, Spain
  (\texttt{fbermejo@dsic.upv.es}).}
\and Blanca Mellado-Pinto\thanks{D.\ Sistemes Inform\`atics i Computaci\'o, Universitat Polit\`ecnica de Val\`encia, Val\`encia, Spain
  (\texttt{bmelpin@dsic.upv.es}).}
\and Jose E. Roman\thanks{D.\ Sistemes Inform\`atics i Computaci\'o, Universitat Polit\`ecnica de Val\`encia, Val\`encia, Spain
  (\texttt{jroman@dsic.upv.es}).}
}

\begin{document}

\maketitle

\begin{abstract}
The non-Hermitian Bethe--Salpeter eigenvalue problem, in the definite case, is a structured eigenproblem, with real eigenvalues coming in pairs $\{\lambda,-\lambda\}$ where the corresponding pair of eigenvectors are closely related, and furthermore the left eigenvectors can be trivially obtained from the right ones. We exploit these properties to devise three variants of structure-preserving Lanczos eigensolvers to compute a subset of eigenvalues (those of either smallest or largest magnitude) together with their corresponding right and left eigenvectors. For this to be effective in real applications, we need to incorporate a thick-restart technique in a way that the overall computation preserves the problem structure. The new methods are validated in an implementation within the SLEPc library using several test matrices, some of them coming from the \yambo materials science code.
\end{abstract}

\section{Introduction}\label{sec:intro}

We are concerned with the eigenvalue problem $Hx=\lambda x$, where the matrix has the form
\begin{equation}
\label{eq:bse1}
H=
\begin{bmatrix}
R & C\\
-C^* & -R^T
\end{bmatrix}\in\mathbb{C}^{2n\times 2n},
\end{equation}
with $R,C\in\mathbb{C}^{n\times n}$ and $R=R^*$, $C=C^T$. We use the $\cdot^T$ and $\cdot^*$ superscripts to denote transposition and complex conjugate transposition, respectively. In addition to eigenvalues $\lambda$ and right eigenvectors $x$, we are also interested in computing left eigenvectors $y$, i.e., those vectors satisfying the relation $y^*H=\lambda y^*$. Furthermore, when a certain condition is met, to be discussed later on, all eigenvalues are real.

Matrix $H$ in~\eqref{eq:bse1} is called the Bethe--Salpeter Hamiltonian matrix as it stems from the Bethe--Salpeter equation~\cite{Salpeter:1951:REB}. This equation is relevant in the analysis of optical absorption and emission processes in solid-state systems by means of first-principles calculation based on Green’s function theory. This is one of the main functionalities offered by the \yambo simulation code~\cite{Marini:2009:YAI,Attaccalite:2011:RAO,Sangalli:2019:MPT}. Starting from a Density Functional Theory (DFT) approximation, \yambo is able to perform additional computations to formulate a Bethe--Salpeter eigenproblem from which to obtain excited-state properties such as energies of quasiparticles. An example of quasiparticle is the so called \emph{exciton}, consisting of a bound state of an electron and a hole; understanding the behaviour of excitons gives valuable information about optical properties of new materials such as two-dimensional extended structures. Other simulation codes with similar goals are BerkeleyGW~\cite{Deslippe:2012:BMP}, MOLGW~\cite{Bruneval:2016:M1M}, and VOTCA-XTP~\cite{Tirimbo:2020:EES}. 
Some of these codes are restricted to a particular case of the BSE problem called the linear response eigenvalue problem, as will be discussed in \cref{sec:bse}. With these tools, it is possible to predict the optical behavior of materials, which is important to design more efficient photovoltaic and light-emitting devices, and in many other applications.

The diagonal blocks of $H$ are associated with the resonant ($R$) and anti-resonant transitions, i.e., from occupied to empty states and from empty to occupied, respectively. In the general case, the anti-resonant terms may be unrelated to $R$, but in most types of analyses both blocks have a certain mutual symmetry relation, and this is what confers special algebraic properties on $H$. The off-diagonal block ($C$) is related to coupling between the two types of transitions. If coupling is small, then these terms can be neglected and the problem can be simplified to a Hermitian eigenproblem for $R$ (the so called Tamm-Dancoff approximation). However, in many interesting applications keeping the coupling is of key importance and hence the problem must be formulated as in~\eqref{eq:bse1}.

For most physical systems, $H$ is a \emph{definite} Bethe--Salpeter Hamiltonian matrix. This will be defined formally in the next section, but it essentially means that if we remove the minus signs in the bottom block-row, the resulting matrix is Hermitian and positive definite. We will assume this property throughout the paper. In that case, all eigenvalues $\lambda$ are real, as will be shown with more detail later on.
However, it is possible to find in the literature cases where $H$ is not definite Bethe-Salpeter, for example in models that include non-negative energy states from the Dirac sea of vacuum~\cite{Ring:2001:TDR}, or because of an instability of the model producing non-physical excitations~\cite{Rangel:2017:ALE}. Benner and Penke~\cite{Benner:2022:EAA} describe a framework for the full eigendecomposition of $H$ that holds even without the definiteness property.

The Bethe--Salpeter equation belongs to a specific class of structured eigenvalue problems. Structured eigenproblems~\cite{Kressner:2005:NMG,Fassbender:2006:SEP} are those whose defining matrices are structured, i.e., their $N^2$ entries depend on less than $N^2$ parameters (with $N=2n$ in our case). Apart from the most obvious Hermitian class, we can find many other classes of structured matrices such as skew-Hermitian, Hamiltonian, unitary, and symplectic, to name a few~\cite{Bunse-Gerstner:1992:CNM}. Preserving the structure can help preserve physically relevant symmetries in the eigenvalues of the matrix and may improve the accuracy and efficiency of the eigensolver. For example, in quadratic eigenvalue problems arising from gyroscopic systems, eigenvalues appear in quadruples $\{\lambda,-\lambda,\bar\lambda,-\bar\lambda\}$, i.e., the spectrum is symmetric with respect to both the real and imaginary axes. This problem can be linearized to a $2\times 2$-block Hamiltonian/skew-Hamiltonian matrix pencil with the same eigenvalues. A structure-preserving eigensolver may give a more accurate answer because it enforces the structure throughout the computation. We can find in the literature many proposals to adapt both direct and iterative methods to particular classes of structured eigenproblems. For instance, in the last decades Krylov methods for Hamiltonian and related eigenproblems have been developed~\cite{Mehrmann:2001:SMC,Watkins:2004:HSL,Benner:2008:SLQ,Benner:2011:HMB}.

Matrix $H$ belongs to the class of complex $J$-symmetric matrices~\cite{Bunse-Gerstner:1992:CNM,Benner:2018:SRC}, where $J=\big[\begin{smallmatrix}0 & I\\ -I & 0\end{smallmatrix}\big]$, with eigenvalues appearing in pairs $\{\lambda,-\lambda\}$. Additionally, it also belongs to the class of pseudo-Hermitian matrices, whose eigenvalues appear in pairs $\{\lambda,\bar{\lambda}\}$. Therefore the eigenvalues of $H$ appear generally as purely real or purely imaginary pairs $\{\lambda,-\lambda\}$ or in complex quadruples $\{\lambda,-\lambda,\bar\lambda,-\bar\lambda\}$, while the definite case enforces a real spectrum $\{\lambda,-\lambda\}$. The eigenvectors also exhibit structural properties in both cases. This will be presented in detail for the definite case in~\cref{sec:problemdef}.

Several authors have recently addressed the Bethe--Salpeter eigenvalue problem from the perspective of structure-preserving eigensolvers. Penke and co-authors~\cite{Penke:2020:HPS,Penke:2022:EAS} have devised effective dense eigensolvers to compute the full eigendecomposition of $H$. On the other hand, Shao and co-authors also focused initially on the computation of the full eigenspectrum with direct methods~\cite{Shao:2016:SPP}, but then turned their attention to iterative methods~\cite{Shao:2018:SPL}, which are able to approximate the optical absorption spectrum without having to compute all eigenvalues. This latter work builds upon previous contributions from the computational physics community~\cite{Gruning:2011:ITL,Haydock:1980:RSS}, providing a more formal view in terms of numerical analysis. The article of Shao \emph{et al.}~\cite{Shao:2018:SPL} is also the starting point of our work, and has been our inspiration to develop the three methods proposed in~\cref{sec:thick}. The key difference between our work and the cited papers is that we are interested in using Lanczos-type methods to explicitly compute a subset of eigenvalues with the corresponding (right and left) eigenvectors, whereas the previous works employ Lanczos recurrences combined with quadrature rules to evaluate the matrix function that is mathematically related to the definition of the optical absorption spectrum. In our case, the eigenvalues and eigenvectors can then be used to reconstruct the optical spectrum, but in addition they provide more information about the spatial distribution of the wave function associated with each excited energy.

Our main contribution is the addition of restart techniques to the Lanczos methods for Bethe--Salpeter. The Lanczos algorithms of Shao \emph{et al.}~\cite{Shao:2018:SPL} and Gr\"uning \emph{et al.}~\cite{Gruning:2011:ITL} rely on a plain Lanczos recurrence. It is well-known that this simple recurrence suffers from loss of orthogonality among the Lanczos vectors, resulting in spurious duplicates of approximate eigenvalues. This phenomenon is benign in the case of approximating the optical absorption spectrum directly. In contrast, when we want to compute eigentriplets accurately, it is necessary to avoid the loss of orthogonality by some kind of reorthogonalization. As a consequence, it is compulsory to keep all the previously computed Lanczos vectors, which increases memory cost in a practical implementation. Restarting techniques aim at alleviating this problem, by bounding the maximum number of vectors to be kept. The thick-restart Lanczos method~\cite{Wu:2000:TLM} is a very effective restart technique that is able to compress the basis of Lanczos vectors while preserving the most relevant spectral information generated so far. It is a particular instance of the more general Krylov--Schur method for non-Hermitian matrices~\cite{Stewart:2001:KAL}.

Our goal is to adapt the thick-restart Lanczos technique to the Bethe--Salpeter eigenproblem, as we have done before in other classes of problems such as the SVD~\cite{Hernandez:2008:REP} or the GSVD~\cite{Alvarruiz:2024:TLB}. The challenge in the case of Bethe--Salpeter is that the restarting step, in addition to preserving the wanted eigenvalues in the compressed basis, must do so in a way that it keeps the problem structure, so that the full computation is still structure-preserving. We will show how to do this in three flavours of Lanczos methods, in \cref{sec:thick}.
In~\cref{sec:unified} we will show how these three different methods can be described using a general framework in terms of Krylov decompositions that explains how structural properties are preserved by the proposed restart mechanisms, and \cref{sec:error} will briefly discuss error bounds.

We have implemented our methods in SLEPc, the Scalable Library for Eigenvalue Problem Computations~\cite{Hernandez:2005:SSF}. \Cref{sec:implem} discusses some implementation details. SLEPc is freely available and widely used in many scientific computing projects around the world. In particular, the \yambo code has the capability to solve the Bethe--Salpeter equation via SLEPc. The development of the new solvers has been done in close collaboration with the \yambo developers. \yambo users will seamlessly benefit from the new solvers' increased efficiency and improved accuracy, being able to pursue more challenging investigations. As examples of promising materials analyzed by \yambo users with the help of SLEPc, we can mention bismuth triiodide~\cite{Cervantes-Villanueva:2024:ELB}, black phosphorus~\cite{Thomen:2024:SSR}, and metallic MXene multilayers~\cite{Kandemir:2024:OPM}. In \cref{sec:results} we illustrate the performance of the new solvers with several test cases, some of them coming from real analyses using \yambo.

\section{Problem definition}\label{sec:problemdef}

In this section we provide more details about the Bethe--Salpeter matrix and its properties, and review previous efforts to apply Lanczos methods to it in a structure-preserving way.

\subsection{The Bethe--Salpeter eigenproblem}\label{sec:bse}

Since $R$ is Hermitian and $C$ is complex symmetric, matrix $H$ in~\eqref{eq:bse1} can also be expressed as
\begin{equation}
\label{eq:bse2}
H=
\begin{bmatrix}
R & C\\
-\overline{C} & -\overline{R}
\end{bmatrix},
\end{equation}
where the bar indicates complex conjugation. In some applications, e.g., systems with real-space inversion symmetry, $R$ and $C$ are real and hence the Bethe--Salpeter matrix is also real, $H=\left[\begin{smallmatrix} R & C \\ -C & -R \end{smallmatrix}\right]$. Shao and coauthors~\cite{Shao:2018:SPL} propose a specific method for this case. This is called the linear response eigenvalue problem (LREP) and has been studied in a series of papers by Bai and Li~\cite{Bai:2012:MPL,Bai:2013:MPL}. Benner and Penke~\cite{Benner:2022:EAA} also analyze the linear response eigenvalue problem, but with complex $R$ and $C$. In our case, we consider only the form~\eqref{eq:bse1} with complex matrices, but we remark that the methods proposed in \cref{sec:thick} can be applied to real matrices without any change.

To derive the algorithms, it will be helpful to factor out the signs in $H$,
\begin{equation}
\label{eq:factored}
H=\Hsign\Hhat,\qquad
\Hsign=\begin{bmatrix}
I & 0\\
0 & -I
\end{bmatrix},
\quad
\Hhat=\begin{bmatrix}
R & C\\
\overline{C} & \overline{R}
\end{bmatrix},
\end{equation}
where $I$ denotes the $n\times n$ identity matrix, $\Hsign$ is a signature matrix and $\Hhat$ is Hermitian. As mentioned in \cref{sec:intro}, we are interested in the case when $\Hhat$ is positive definite, in which case we say $H$ is a definite Bethe--Salpeter matrix. Throughout the paper, we will assume that this property holds, which is generally true in most applications.

The Bethe--Salpeter matrix $H$ belongs to the slightly more general class of complex $J$-symmetric matrices, that can be written as
\begin{equation}
\label{eq:jsymmetric}
H_C=
\begin{bmatrix}
R & C\\
D & -R^T
\end{bmatrix}\in\mathbb{C}^{2n\times 2n},
\quad R,C=C^T,D=D^T\in\mathbb{C}^{n\times n}.
\end{equation}
These matrices satisfy $JH_C=(JH_C)^T$, where $J=\big[\begin{smallmatrix}0 & I\\ -I & 0\end{smallmatrix}\big]$.
The eigenvalues of complex $J$-symmetric matrices display a symmetry~\cite{Bunse-Gerstner:1992:CNM,Benner:2018:SRC}: they appear in pairs $\{\lambda,-\lambda\}$. Also, if $x$ is the right eigenvector for $\lambda$, $H_Cx=\lambda x$, then $J\bar{x}$ is the left eigenvector for $-\lambda$,
\begin{equation}
\label{eq:leftevec}
(J\bar{x})^*H_C=-\lambda(J\bar{x})^*.
\end{equation}

The Bethe--Salpeter matrix $H$ inherits these properties, but has some additional structure. Since $H=S\Hhat$, with $\Hhat$ Hermitian positive definite, it can be transformed by similarity to the Hermitian matrix $\Hhat^{1/2}S\Hhat^{1/2}$, thus $H$ is diagonalizable and has real spectrum.
Due to this fact, a specific version of the property in~\eqref{eq:leftevec} can be proved~\cite{Shao:2016:SPP}:

\begin{theorem}[Shao \emph{et al.}~\cite{Shao:2016:SPP}]\label{thm:shao}
Let H be of the form~\eqref{eq:bse1} satisfying that $\Hhat$ in~\eqref{eq:factored} is positive definite. Then there exist $X_1$, $X_2\in\mathbb{C}^{n\times n}$ and positive numbers $\lambda_1,\lambda_2,\dots,\lambda_n\in\mathbb{R}$ such that
\begin{equation}
\label{eq:eigendecomposition}
HX=X
\begin{bmatrix}
\Lambda_+ & 0 \\
0 & -\Lambda_+
\end{bmatrix},\qquad
Y^*H=
\begin{bmatrix}
\Lambda_+ & 0 \\
0 & -\Lambda_+
\end{bmatrix}Y^*,\qquad
Y^*X=I_{2n},
\end{equation}
where
\begin{equation}
\label{eq:spliteigenvectors}
X=
\begin{bmatrix}
X_1 & \overline{X}_2 \\
X_2 & \overline{X}_1
\end{bmatrix},\qquad
Y=
\begin{bmatrix}
X_1 & -\overline{X}_2 \\
-X_2 & \overline{X}_1
\end{bmatrix},
\end{equation}
and $\Lambda_+=\operatorname{diag}\{\lambda_1,\dots,\lambda_n\}$.
\end{theorem}

The previous theorem provides a complete description of the eigendecomposition of $H$, showing that there is considerable redundancy. Efficient eigensolvers should exploit this fact to compute the solution faster.

For large-scale problems, computing the full eigendecomposition is too costly. For an accurate enough approximation of the optical absorption spectrum, it is sufficient to compute a few eigenvalues and their corresponding eigenvectors. That is why we focus on iterative eigensolvers, particularly Lanczos methods that can provide good approximations to exterior eigenvalues in relatively few iterations. Here we remark that matrices coming from Green's function theory such as those in the \yambo code are dense, while Lanczos is most appropriate in the case of sparse matrices, or more generally when matrix-vector products are cheap. Still, this approach can be competitive with respect to full diagonalization, provided that the number of wanted eigenvalues is not too large. In this paper, we assume that the user is interested in a few hundreds of eigenvalues at most.

The eigenvalues of interest in applications are those of smallest magnitude, i.e., the smallest energies. Due to the $\{\lambda,-\lambda\}$ symmetry, those eigenvalues lie in the middle of the spectrum, which is a bad scenario for the convergence of Lanczos methods. However, the variants of Lanczos presented in~\cref{sec:thick} work around this issue by some kind of \emph{spectrum folding} approach, where negative eigenvalues become positive. Essentially, the proposed methods compute only $\lambda\in\Lambda_+$ with corresponding eigenvectors, and then use \cref{thm:shao} to inexpensively reconstruct the spectral elements for $-\lambda$, including left eigenvectors which are also required by the application. Hence, techniques such as the shift-and-invert spectral transformation for computing interior eigenvalues are not usually necessary, unless eigenvalues are clustered together close to zero. Still, in \cref{sec:implem} we discuss how this technique can be implemented in our Bethe--Salpeter solvers.

\subsection{Possible solution approaches and related work}\label{sec:related}

The most naive solution approach is to use a general non-Hermitian Krylov solver, configuring it to search for smallest magnitude eigenvalues. It must be capable of computing left eigenvectors as well. For this, SLEPc provides an implementation of two-sided Krylov--Schur~\cite{Zwaan:2017:KRT}. This was the strategy used by the \yambo code prior to the work presented in this paper. A better approach would be to use plain Krylov--Schur to obtain right eigenvectors only, and then apply \cref{thm:shao} to obtain left eigenvectors. This will be the baseline for the performance comparisons in \cref{sec:results}. However, this approach is neither enforcing the structure of the problem, for example the realness of the eigenvalues or the symmetry of the spectrum, nor taking it into account to improve efficiency.

The Bethe--Salpeter eigenvalue problem can be converted to a real Hamiltonian eigenvalue problem~\cite[Thm.~2]{Shao:2016:SPP}, and this connection has been exploited by some authors to develop full diagonalization eigensolvers~\cite{Shao:2016:SPP,Penke:2020:HPS}. Similarly, we could adapt the structure-preserving Krylov methods for Hamiltonian eigenvalue problems mentioned in \cref{sec:intro}. However, this would exploit only part of the structure, and disregard the richer structure given by \cref{thm:shao}. Although we do not consider this route here, we remark that it would be a potentially good approach in case of non-definite problems, i.e., when $\Hhat$ is not definite, in which case the methods of \cref{sec:thick} cannot be applied.

In view of~\eqref{eq:factored}, we could consider a methodology for the product eigenvalue problem. The details for periodic Krylov--Schur have been worked out by Kressner~\cite{Kressner:2006:PKA}. If we applied this method to~\eqref{eq:factored}, we would obtain a periodic Krylov decomposition with two blocks, where the projected problem would contain an upper Hessenberg and an upper triangular block, to which it is possible to apply a Krylov--Schur-type restart. Interestingly, the Lanczos decompositions that we will describe in detail in \cref{sec:thick} resemble this, but with tridiagonal and diagonal blocks since we additionally exploit the symmetries underlying the problem structure.

To understand these symmetries, we have to interpret~\eqref{eq:factored} as a generalized eigenvalue problem. If we consider the pencil $(\Hsign,\Hhat)$, we have a Hermitian definite generalized eigenproblem where we can apply a $B$-Lanczos recurrence with the inner product induced by $\Hhat$. The Bethe--Salpeter matrix $H$ is self-adjoint with respect to this inner product, and Lanczos will produce a real symmetric tridiagonal matrix containing Ritz approximations of the eigenvalues of $H$. An alternative would be to consider the matrix pair $(\Hhat,\Hsign)$. It can be shown that $H$ is also self-adjoint with respect to the indefinite inner product induced by $\Hsign$, so an analogous approach could be applied using a pseudo-Lanczos recurrence. We note that SLEPc provides an implementation of pseudo-Lanczos that can be useful in some settings~\cite{Campos:2016:RQM}, but it is better to avoid it since numerical stability cannot be guaranteed. The approach with definite inner product is to be preferred, even though inner products for $\Hhat$ are much more expensive than those for $\Hsign$. The good news is that the methods presented below get around the high cost of $\Hhat$-inner products and never need to compute them explicitly.

Shao and coauthors~\cite{Shao:2018:SPL} propose a structure-preserving Lanczos process using $\Hhat$-inner products. There is a connection between the quantities computed by this method and a Lanczos decomposition for $H^2$. This implicit squaring of the matrix has the effect of folding the spectrum around zero so that all eigenvalues become positive. This method is described in detail in \cref{sec:shao}. The last part of~\cite{Shao:2018:SPL} analyzes the connection of their method with other Lanczos methods. One of them is due to Gr\"uning \emph{et al.}~\cite{Gruning:2011:ITL} and we will discuss it in \cref{sec:gruning}. Another variation is also analyzed in connection with symplectic Lanczos, where the projected eigenproblem of the Lanczos factorization has again a Bethe--Salpeter structure, see the details in \cref{sec:projectedbse}. As mentioned in \cref{sec:intro}, all these methods were originally intended for computing the absorption spectrum directly via quadrature rules. In contrast, our presentation of the methods in \cref{sec:thick} focuses on computation of eigenvalues and eigenvectors, with an effective restart.

We conclude this section by mentioning an additional strategy for solving the Bethe--Salpeter eigenproblem. It is possible to build an equivalent product eigenvalue problem defined by two real matrices of order $2n$~\cite{Zhang:2024:CLO}. A similar product eigenvalue problem is often used as a proxy for solving the LREP, so this opens the door to using, in the context of the Bethe--Salpeter eigenproblem, methods that have been successfully employed in the LREP, such as several variants of LOBPCG~\cite{Zhang:2024:CLO}.

\section{Thick-restart Lanczos}\label{sec:thick}

The three methods discussed in this section are mathematically equivalent, since they are based on equivalent Lanczos recurrences, only differing in the way the different quantities are computed. Computational results in \cref{sec:results} will show that all of them provide good accuracy and performance.

Throughout this section we will use the following notation. The $k\times k$ identity matrix is denoted by $I_k$ and the $j$th column of the identity is denoted by $e_j$, where the length of the vector is implicitly defined by the context. For matrices representing basis vectors such as $V_k$, the subindex denotes the number of columns, while in square matrices such as the tridiagonal $T_k$, the subindex refers to dimension in both rows and columns. The non-standard inner product induced by a Hermitian positive-definite matrix $B$ is denoted by $(x,y)_B=y^*Bx$, and similarly the corresponding norm is $\|x\|_B=\sqrt{x^*Bx}$. The imaginary unit will be denoted by i (in roman) and $\text{Re}(\cdot)$ and $\text{Im}(\cdot)$ will respectively denote the real and imaginary parts of a complex entity.

\subsection{Method based on Shao \emph{et al}}\label{sec:shao}
\def\ShaoW{G}
\subsubsection{Lanczos recurrence}
The Lanczos procedure proposed in \cite{Shao:2018:SPL} builds the Lanczos-type relation
\begin{equation}
H\begin{bmatrix}
U_k & V_k\\
\mybar{U}_k & -\mybar{V}_k
\end{bmatrix}
=
\begin{bmatrix}
U_k & V_k\\
\mybar{U}_k & -\mybar{V}_k
\end{bmatrix}
\begin{bmatrix}
0 & T_k\\I_k&0
\end{bmatrix}
+\beta_k\begin{bmatrix}
u_{k+1}\\\bar{u}_{k+1}
\end{bmatrix}
e_{2k}^*,
\label{eq:shao-lanczos}
\end{equation}
where $U_k,V_k\in\mathbb{C}^{n\times k}$, $U_k=[u_1,\dots,u_k],V_k=[v_1,\dots,v_k]$, and $T_k$ is the real symmetric tridiagonal positive definite matrix
\[
T_k=\text{tridiag}
\left\lbrace
\begingroup 
\setlength\arraycolsep{2pt}
\begin{matrix}
        &\beta_1&      &\cdots&      &\beta_{k-1}& \\
\alpha_1&       &\cdots&      &\cdots&           &\alpha_{k}\\
        &\beta_1&      &\cdots&      &\beta_{k-1}&
\end{matrix}
\endgroup
\right\rbrace.
\]

Lanczos vectors satisfy the orthogonality condition
\begin{equation}
\begin{bmatrix}
V_k      & U_k\\
\mybar{V}_k&-\mybar{U}_k
\end{bmatrix}^*
\begin{bmatrix}
U_k & V_k\\
\mybar{U}_k & -\mybar{V}_k
\end{bmatrix}
=
2I_{2k}.
\label{eq:shao-ortog}
\end{equation}

Assuming $u_{j-1},u_j,v_j,\beta_{j-1}$, are available for a given $j\ge1$, with $u_0=0$, new vectors are computed by first generating
\begin{equation}\label{eq:shao-recurrence}
\begin{aligned}
\tilde{u}_{j+1}
&=
Rv_j-C\bar{v}_j
-\beta_{j-1}u_{j-1}
-\alpha_j u_{j},
\\
\tilde{v}_{j+1}
&=
R\tilde{u}_{j+1}+C\bar{\tilde{u}}_{j+1},
\end{aligned}
\end{equation}
where $\alpha_j=\text{Re}(v_j^*(Rv_j-C\bar{v}_j))$, and then normalizing to obtain
\[
u_{j+1}=\tilde{u}_{j+1}/\beta_j,\quad v_{j+1}=\tilde{v}_{j+1}/\beta_j,
\qquad
\text{with }\beta_j=\sqrt{\text{Re}(\tilde{u}_{j+1}^*\tilde{v}_{j+1})}.
\]

Initial vectors can be generated from an arbitrary $\tilde{u}_1$, by computing $\tilde{v}_{1}=R\tilde{u}_{1}+C\bar{\tilde{u}}_{1}$ and normalizing
\[
u_{1}=\tilde{u}_{1}/\beta_0,\quad v_{1}=\tilde{v}_{1}/\beta_0,
\qquad
\text{with }\beta_0=\sqrt{\text{Re}(\tilde{u}_{1}^*\tilde{v}_{1})}.
\]

As shown in \cite{Shao:2018:SPL}, matrices $U_k$ and $T_k$ can equivalently be obtained by applying a $k$-step Lanczos procedure to $H^2$ with the $\Hhat$-inner product, resulting in
\[
H^2\begin{bmatrix}U_k\\\mybar{U}_k\end{bmatrix}
=
\begin{bmatrix}U_k\\\mybar{U}_k\end{bmatrix}T_k
+\beta_k \begin{bmatrix}u_{k+1}\\\bar{u}_{k+1}\end{bmatrix} e_k^*,
\]
with the orthogonality condition corresponding to \eqref{eq:shao-ortog} being
\[
\begin{bmatrix}
U_k\\
\mybar{U}_k 
\end{bmatrix}^*
\Hhat
\begin{bmatrix}
U_k\\
\mybar{U}_k
\end{bmatrix}
=
2I_{k}.
\]

\subsubsection{Reorthogonalization}\label{sec:shao-reorth}
As mentioned in \cref{sec:intro}, loss of orthogonality in the Lanczos recurrence to build~\eqref{eq:shao-lanczos} is not a problem when using it to compute the absorption spectrum, as in \cite{Shao:2018:SPL}, but must be avoided in our case, because our goal is to compute the eigenvalues and eigenvectors.
Our approach is to do full reorthogonalization at each Lanczos step, see \cref{sec:implem} for additional discussion.
We show next how to implement reorthogonalization in Shao's method.

Suppose orthogonality of vector $\tilde{u}_{j+1}$ is not perfect, so that we have
\[
\beta_j\begin{bmatrix}u_{j+1}\\\bar{u}_{j+1}\end{bmatrix}
=
\begin{bmatrix}\tilde{u}_{j+1}\\\bar{\tilde{u}}_{j+1}\end{bmatrix}
-\begin{bmatrix}U_j&V_j\\\mybar{U}_{j}&-\mybar{V}_{j}\end{bmatrix}
\begin{bmatrix}c\\d\end{bmatrix},
\]
with $c,d$ not zero.
Premultiplying by
$\left[\begin{smallmatrix}
V_j      & U_j\\
\bar{V}_j&-\bar{U}_j
\end{smallmatrix}\right]^*
$ and taking into account the orthogonality condition \eqref{eq:shao-ortog}, we get
\[
\begin{bmatrix}
c\\d
\end{bmatrix}
=
\frac{1}{2}
\begin{bmatrix}
V_j      & U_j\\
\mybar{V}_j&-\mybar{U}_j
\end{bmatrix}^*
\begin{bmatrix}\tilde{u}_{j+1}\\\bar{\tilde{u}}_{j+1}\end{bmatrix},
\]
which leads to
\[
\begin{aligned}
c&=\text{Re}(V_j^*\tilde{u}_{j+1}),\\
d&=\text{Im}(U_j^*\tilde{u}_{j+1})\text{i}.
\end{aligned}
\]

After computing $c,d$, a reorthogonalized vector $\check{u}_{j+1}$ is obtained as
\[
\check{u}_{j+1}=\beta_ju_{j+1}=\tilde{u}_{j+1}- U_j c - V_j d.
\]

The resulting Lanczos process is described by means of \cref{alg:shao-reorth}.

\begin{algorithm}
  \caption{Shao's Lanczos with reorthogonalization}
  \label{alg:shao-reorth}
  \begin{algorithmic}[1]
    \REQUIRE A definite $2n\times 2n$ BSE matrix $H$ with structure \eqref{eq:bse1}; initial vector $\tilde{u}_1$ of size $n$; number of Lanczos steps $k$.
    \ENSURE $U_k,V_k,T_k$ making up Lanczos-type decomposition \eqref{eq:shao-lanczos}.
    \STATE $\tilde{v}_{1}=R\tilde{u}_{1}+C\bar{\tilde{u}}_{1}$
    \STATE $\beta_0=\sqrt{\text{Re}(\tilde{u}_{1}^*\tilde{v}_{1})},\quad u_0=0$
    \STATE $u_{1}=\tilde{u}_{1}/\beta_0,\quad v_{1}=\tilde{v}_{1}/\beta_0$
    \FOR{$j=1,2,\dots,k$}
    \STATE $x=Rv_j-C\bar{v}_j$
    \STATE $\tilde{\alpha}_j=\text{Re}(v_j^*x)$
    \STATE $\tilde{u}_{j+1}=x-\beta_{j-1}u_{j-1}-\tilde{\alpha}_j u_{j}$
    \STATE $c=\text{Re}(V_j^* \tilde{u}_{j+1}),\quad d=\text{Im}(U_j^* \tilde{u}_{j+1})\text{i}$
    \STATE $\alpha_j = \tilde{\alpha}_j+c_j$
    \STATE $\check{u}_{j+1}=\tilde{u}_{j+1}-U_j c-V_j d$
    \STATE $\check{v}_{j+1}=R\check{u}_{j+1}+C\bar{\check{u}}_{j+1}$
    \STATE $\beta_j=\sqrt{\text{Re}(\check{u}_{j+1}^*\check{v}_{j+1})}$
    \STATE $u_{j+1}=\check{u}_{j+1}/\beta_j,\quad v_{j+1}=\check{v}_{j+1}/\beta_j$
    \ENDFOR
  \end{algorithmic}
\end{algorithm}

\subsubsection{Restart}

Let $T_k=QD_kQ^*$, with $Q$ orthogonal, be the eigendecomposition of $T_k$. Then,
\[
\begin{bmatrix}
0  &T_k\\
I_k & 0
\end{bmatrix}
=
\begin{bmatrix}
Q & 0\\
0 &Q
\end{bmatrix}
\begin{bmatrix}
0  &D_k\\
I_k & 0
\end{bmatrix}
\begin{bmatrix}
Q^* & 0\\
0   &Q^*\textit{}
\end{bmatrix}.
\]

Substituting in \eqref{eq:shao-lanczos} and multiplying on the right by $\left[\begin{smallmatrix}Q & 0\\0 & Q\end{smallmatrix}\right]$, we get
\begin{equation}
H
\begin{bmatrix}
U_k Q         & V_k Q\\
\mybar{U}_k Q & -\mybar{V}_k Q
\end{bmatrix}
=
\begin{bmatrix}
U_k Q         & V_k Q\\
\mybar{U}_k Q & -\mybar{V}_k Q
\end{bmatrix}
\begin{bmatrix}
0  &D_k\\
I_k & 0
\end{bmatrix}
+\beta_k\begin{bmatrix}
u_{k+1}\\\bar{u}_{k+1}
\end{bmatrix}
\begin{bmatrix}
0 &e_{k}^*Q
\end{bmatrix},
\label{eq:shao-lanczos-updated}
\end{equation}
or, defining $\hat{U}_k:=U_kQ$, $\hat{V}_k:=V_kQ$ and $b:=\beta_kQ^*e_{k}$,
\begin{equation}
H
\begin{bmatrix}
\hat{U}_k            & \hat{V}_k\\
\mybar{\hat{U}}_k & -\mybar{\hat{V}}_k
\end{bmatrix}
=
\begin{bmatrix}
\hat{U}_k            & \hat{V}_k\\
\mybar{\hat{U}}_k & -\mybar{\hat{V}}_k
\end{bmatrix}
\begin{bmatrix}
0  &D_k\\
I_k & 0
\end{bmatrix}
+\begin{bmatrix}
u_{k+1}\\\bar{u}_{k+1}
\end{bmatrix}
\begin{bmatrix}
0  &b^*
\end{bmatrix},
\label{eq:shao-lanczos-trunc}
\end{equation}

Since $D_k$ is diagonal, the previous decomposition can be truncated so that each side of the equality contains only $2r$ columns, with $r<k$.
Before truncating, we permute the eigendecomposition of $T_k$ so that the leading block of $D_k$ contains the wanted eigenvalues, i.e., the smallest ones.

After truncation, the decomposition \eqref{eq:shao-lanczos-trunc} can be extended again to $2k$ columns, by running \cref{alg:shao-reorth-reset}, which is a modified version of \cref{alg:shao-reorth} that receives vectors $[\hat{u}_1,\dots,\hat{u}_r,u_{k+1}]$ and $[\hat{v}_1,\dots,\hat{v}_r,v_{k+1}]$ as the columns of $U_{r+1}$ and $V_{r+1}$, respectively, and the vector $b=[b_1,\dots,b_r]^T$. Note that vector $\tilde{u}_{r+2}$, generated in the first iteration of the algorithm, must be orthogonalized against all the columns of $U_{r}$, using $b_1,\dots,b_r$ as the orthogonalization coefficients.
\Cref{alg:shao-reorth-reset} is a generalization of \cref{alg:shao-reorth}, so it can also be used for the Lanczos process before the first restart, setting $r=0$ and $b=[\,]$.

\begin{algorithm}
  \caption{Shao's Lanczos with reorthogonalization and accommodating restart}
  \label{alg:shao-reorth-reset}
  \begin{algorithmic}[1]
    \REQUIRE A definite $2n\times 2n$ BSE matrix $H$ with structure \eqref{eq:bse1}; initial vectors in columns of $U_{r+1}$ and $V_{r+1}$; vector $b$; number of Lanczos steps $k$.
    \ENSURE $U_k,V_k,T_k$ extending decomposition \eqref{eq:shao-lanczos-trunc} to $2k$ columns.
    \FOR{$j=r+1,r+2,\dots,k$}
    \STATE $x=Rv_j-C\bar{v}_j$
    \STATE $\tilde{\alpha}_j=\text{Re}(v_j^*x)$ \label{alg:shao-reorth-reset:line3}
    \IF{$j=r+1$}
    \STATE $\tilde{u}_{j+1}=x-U_{r}b-\tilde{\alpha}_j u_{j}$
    \ELSE
    \STATE $\tilde{u}_{j+1}=x-\beta_{j-1}u_{j-1}-\tilde{\alpha}_j u_{j}$
    \ENDIF \label{alg:shao-reorth-reset:line8}
    \STATE $c=\text{Re}(V_j^* \tilde{u}_{j+1}),\quad d=\text{Im}(U_j^* \tilde{u}_{j+1})\text{i}$ \label{alg:shao-reorth-reset:line9}
    \STATE $\alpha_j = \tilde{\alpha}_j+c_j$
    \STATE $\check{u}_{j+1}=\tilde{u}_{j+1}-U_j c-V_j d$ \label{alg:shao-reorth-reset:line11}
    \STATE $\check{v}_{j+1}=R\check{u}_{j+1}+C\bar{\check{u}}_{j+1}$
    \STATE $\beta_j=\sqrt{\text{Re}(\check{u}_{j+1}^*\check{v}_{j+1})}$
    \STATE $u_{j+1}=\check{u}_{j+1}/\beta_j,\quad v_{j+1}=\check{v}_{j+1}/\beta_j$
    \ENDFOR
  \end{algorithmic}
\end{algorithm}

This process is repeated until the desired number of Ritz pairs have converged. We will see below that the entries of $b$ are related to residual norms of each Ritz pair.

\subsubsection{Obtaining approximate eigentriplets}
\label{sec:eigentriplets}

It remains to see how to obtain the approximate eigenvalues and eigenvectors from the Krylov decomposition \eqref{eq:shao-lanczos-trunc}.
Applying a perfect shuffle symmetric permutation
$\begin{bmatrix}e_1,e_{k+1},e_2,e_{k+2},\dots\end{bmatrix}$,
the projected matrix
\[
\begin{bmatrix}
0  &D_{k}\\
I_k & 0
\end{bmatrix}
\]
gets transformed into a block diagonal matrix with $2\times 2$ diagonal blocks of the form
$
\left[
\begin{smallmatrix}
 0 &d_i\\
 1 &0
\end{smallmatrix}
\right]
$
with $d_i>0$. Each such block admits the eigendecomposition
\[
\begin{bmatrix}
 0 &d_i\\
 1 &0
\end{bmatrix}
\begin{bmatrix}
\sqrt{d_i} & -\sqrt{d_i}\\ 
         1 & 1
\end{bmatrix}
=
\begin{bmatrix}
\sqrt{d_i} & -\sqrt{d_i}\\ 
         1 & 1
\end{bmatrix}
\begin{bmatrix}
\sqrt{d_i}& 0\\
0         &-\sqrt{d_i}
\end{bmatrix}.
\]

Thus, the projected matrix admits the eigendecomposition
\[
\begin{bmatrix}
0  &D_{k}\\
I_k & 0
\end{bmatrix}
\ShaoW
=
\ShaoW
\begin{bmatrix}
D_k^\frac{1}{2}& 0\\
0       &-D_k^\frac{1}{2}
\end{bmatrix},
\]
where
\[
\ShaoW:=\begin{bmatrix}
D_k^\frac{1}{2}&D_k^\frac{1}{2}\\
    I_k  &   -I_k
\end{bmatrix}.
\]

Substituting in \eqref{eq:shao-lanczos-trunc} and multiplying on the right by $\ShaoW$, we get
\begin{equation}
H
\begin{bmatrix}
\tilde{X}_1 & \mybar{\tilde{X}}_2 \\
\tilde{X}_2 & \mybar{\tilde{X}}_1
\end{bmatrix}
=
\begin{bmatrix}
\tilde{X}_1 & \mybar{\tilde{X}}_2 \\
\tilde{X}_2 & \mybar{\tilde{X}}_1
\end{bmatrix}
\begin{bmatrix}
\tilde{\Lambda}_+& 0\\
0       &-\tilde{\Lambda}_+
\end{bmatrix}
+\begin{bmatrix}
u_{k+1}\\\bar{u}_{k+1}
\end{bmatrix}
\begin{bmatrix}
b^*  &-b^*
\end{bmatrix},
\label{eq-aprox-vecprop}
\end{equation}
where we use the definitions
\[
\tilde{X}=
\begin{bmatrix}
\tilde{X}_1 & \mybar{\tilde{X}}_2 \\
\tilde{X}_2 & \mybar{\tilde{X}}_1
\end{bmatrix}
:=
\begin{bmatrix}
\hat{U}_k            &\hat{V}_k\\
\mybar{\hat{U}}_k &-\mybar{\hat{V}}_k
\end{bmatrix}
\ShaoW
,\qquad
\tilde{\Lambda}_+:=D_k^\frac{1}{2}.
\]

Clearly, $\tilde{\Lambda}_+$ and $\tilde{X}$ contain approximations to a subset of eigenvalues of $H$ and the corresponding right eigenvectors, respectively. Approximations to left eigenvectors $\tilde{Y}$ can be obtained trivially from $\tilde{X}$ taking into account \eqref{eq:spliteigenvectors}.

Considering column $i\leq k$ of the relation~\eqref{eq-aprox-vecprop}, we can see that the residual norm of the approximate eigenpair $(\tilde\lambda_i,\tilde{x}_i)$ can be obtained as $\|H\tilde{x}_i-\tilde\lambda_i\tilde{x}_i\|_2=\rho|b_i|$, where $\rho=\|[u_{k+1}^*,\bar{u}_{k+1}^*]^*\|_2$, and we have exactly the same residual norm for the corresponding negative eigenpair in column $i+k$. We consider that an eigenpair is converged when $|b_i|<\mathtt{tol}$ for a given tolerance \texttt{tol}, disregarding $\rho$. More precisely, we check for convergence in order, i.e., if the number of wanted eigenvalues is $n_\mathrm{ev}$, the computation will stop whenever the first $n_\mathrm{ev}/2$ residual norms are below the tolerance. Note that the method returns $n_\mathrm{ev}/2$ positive eigenvalues together with the same amount of negative ones.

The previous paragraph considers the absolute error of the approximate eigenpair. An alternative is to use the relative residual estimates $|b_i|/|\tilde\lambda_i|$, where we are assuming that the approximate eigenvector has been normalized so that $\|\tilde{x}_i\|_2=1$ (see \cref{sec:implem} for a discussion about normalization of eigenvectors). In our implementation, the default is to use a relative convergence criterion, but the user can easily switch to an absolute one if preferred.
\revisedtext{}{The implementation also uses \emph{locking} by default, so that converged eigenpairs are frozen in subsequent restarts, although the basis vectors of these converged eigeinpairs still take part in orthogonalization processes.}

The overall restarted method is illustrated in~\cref{alg:shao-restarted}.

\begin{algorithm}
  \caption{Restarted Shao's Lanczos method}
  \label{alg:shao-restarted}
  \begin{algorithmic}[1]
    \REQUIRE A definite $2n\times 2n$ BSE matrix $H$ with structure \eqref{eq:bse1}; initial vector $\tilde{u}_1$ of size $n$; maximum number of Lanczos steps $k$, restart parameter $r$, number of wanted eigenvalues $n_\mathrm{ev}$, tolerance \texttt{tol}.
    \ENSURE Approximate eigenpairs $(\tilde\lambda_i,\tilde{x}_i)$, $i=1,\dots,n_\mathrm{ev}/2$
    \REPEAT
    \STATE Build Lanczos relation of order $k$ with \cref{alg:shao-reorth-reset}
    \STATE Compute sorted eigendecomposition $T_k=QD_kQ^*$
    \STATE Update Lanczos relation as in \eqref{eq:shao-lanczos-updated}
    \STATE Truncate Lanczos relation to order $r$ as in \eqref{eq:shao-lanczos-trunc}
    \STATE Set $n_\mathrm{conv}$ such that $|b_i|<\mathtt{tol}\cdot|\tilde\lambda_i|,\;i=1,\dots,n_\mathrm{conv}$
    \UNTIL{$n_\mathrm{conv}\geq n_\mathrm{ev}/2$}
    \STATE Compute eigendecomposition of $\big[\begin{smallmatrix}0 &D_{k}\\ I_k &0\end{smallmatrix}\big]$, update Lanczos relation as in \eqref{eq-aprox-vecprop}
  \end{algorithmic}
\end{algorithm}

\subsection{Method based on Gr\"uning \emph{et al}}\label{sec:gruning}
\def\GruningU{M}
\def\Gruningu{m}
\def\GruningW{N}
\def\Gruningw{n}
\def\GruningV{V}
\def\Gruningv{v}
\def\GrHU{M'}
\def\GrHW{N'}
\def\GrHu{m'}
\def\GrHw{n'}
\def\GruningLowerBidiagonal{L}
\def\GruningQ{Q}
\def\GruningZ{P}
\def\GruningLambda{\Sigma}
\def\GruningP{\hat{G}}
\def\Grbeta{\hat{\beta}}
\subsubsection{Lanczos recurrence}

The method proposed in~\cite{Gruning:2011:ITL} is a Lanczos iteration with an inner product with $\Hhat$, although it is reformulated in terms of a much cheaper inner product with the diagonal matrix $S$, taking into account that $\Hhat=\Hsign H=H^*\Hsign$.
In particular, for a vector basis $V_j=[v_1,v_2,\dots,v_j]$, the orthogonalization of a new vector $\tilde{\Gruningv}_{j+1} = Hv_{j}$ against the previously computed $V_{j}$ requires computing the inner product $(\tilde{\Gruningv}_{j+1},\GruningV_{j})_{\Hhat}=(H\Gruningv_{j},H\GruningV_{j})_\Hsign$. In the Lanczos iteration, $H\GruningV_{j}$ is not computed for every inner product. Instead, it is stored in memory and a new vector is appended at each iteration as it is computed.

Further optimization is possible using an initial vector with the structure
$\Gruningv_1=\begin{bmatrix}\Gruningu_1\\\bar{\Gruningu}_1\end{bmatrix}$,
as the (non-restarted) Lanczos algorithm produces a sequence of vectors with the form
\[
\GruningV_{2k}=
\begin{bmatrix}
     \Gruningu_1 &      \Gruningw_1 &      \Gruningu_2 &      \Gruningw_2 &\cdots& \Gruningu_{k} & \Gruningw_{k}    \\
\bar{\Gruningu}_1&-\bar{\Gruningw}_1& \bar{\Gruningu}_2&-\bar{\Gruningw}_2&\cdots&\bar{\Gruningu}_{k}& -\bar{\Gruningw}_{k}
\end{bmatrix},
\]
and a tridiagonal matrix $T_{2k}$ with zeros on the diagonal,
\[
T_{2k}=\text{tridiag}
\left\lbrace
\begingroup 
\setlength\arraycolsep{2pt}
\begin{array}{ccccccccccccc}
  &\Grbeta_1&   &\Grbeta_{k+1}&   &\Grbeta_2&        &\cdots&        &\Grbeta_{2k-1}&   &\Grbeta_k &   \\
0 &       & 0 &           & 0 &       & \cdots &      & \cdots &            & 0 &        & 0 \\
  &\Grbeta_1&   &\Grbeta_{k+1}&   &\Grbeta_2&        &\cdots&        &\Grbeta_{2k-1}&   &\Grbeta_k & 
\end{array}
\endgroup
\right\rbrace,
\]
satisfying the relation
\begin{equation}
H\GruningV_{2k}=\GruningV_{2k}T_{2k}+\Grbeta_{2k}\Gruningv_{2k+1}e^*_{2k}.
\label{eq:gruning-before-permutation}
\end{equation}

Postmultiplying the factorization \eqref{eq:gruning-before-permutation} by an odd-even permutation matrix 
gives the resulting Lanczos-type relationship
\begin{equation}
H\begin{bmatrix}
\GruningU_k & \GruningW_k\\
\mybar{\GruningU}_k & -\mybar{\GruningW}_k
\end{bmatrix}
=
\begin{bmatrix}
\GruningU_k & \GruningW_k\\
\mybar{\GruningU}_k & -\mybar{\GruningW}_k
\end{bmatrix}
\begin{bmatrix}
0 & \GruningLowerBidiagonal_k \\ \GruningLowerBidiagonal^*_k & 0
\end{bmatrix}
+\Grbeta_{2k}\begin{bmatrix}
\Gruningu_{k+1}\\\mybar{\Gruningu}_{k+1}
\end{bmatrix}
e_{2k}^*,
\label{eq:gruning-relationship}
\end{equation}
where $\GruningU_k,\GruningW_k\in\mathbb{C}^{n\times k}$, $\GruningU_k=[\Gruningu_1,\dots,\Gruningu_k],\GruningW_k=[\Gruningw_1,\dots,\Gruningw_k]$, and $\GruningLowerBidiagonal_k$ is the real lower bidiagonal positive matrix
\[
\GruningLowerBidiagonal_k=
\begin{bmatrix}
\Grbeta_1 & & & &\\
\Grbeta_{k+1} & \Grbeta_2 & & &\\
& \Grbeta_{k+2} & \Grbeta_3 & & &\\
& & \ddots & \ddots & \\
& & & \Grbeta_{2k-1} & \Grbeta_k \\
\end{bmatrix}.
\]
 
Lanczos vectors satisfy the orthogonality condition
\begin{equation}
\begin{bmatrix}
\GruningU_k & \GruningW_k\\
\mybar{\GruningU}_k & -\mybar{\GruningW}_k
\end{bmatrix}^*
\Hhat
\begin{bmatrix}
\GruningU_k & \GruningW_k\\
\mybar{\GruningU}_k & -\mybar{\GruningW}_k
\end{bmatrix}
=
I_{2k},
\label{eq:gruning-orthog-Hhat}
\end{equation}
which can also be expressed as
\begin{equation}
\begin{bmatrix}
\GrHU_k & \GrHW_k\\
-\mybar{\GrHU}_k & \mybar{\GrHW}_k
\end{bmatrix}^*
\Hsign
\begin{bmatrix}
\GruningU_k & \GruningW_k\\
\mybar{\GruningU}_k & -\mybar{\GruningW}_k
\end{bmatrix}
=
I_{2k},
\label{eq:gruning-orthog}
\end{equation}
where $\GrHU_k=R\GruningU_{k}+C\mybar{\GruningU}_{k}$ and $\GrHW_k=R\GruningW_{k}-C\mybar{\GruningW}_{k}$.

Relation~\eqref{eq:gruning-relationship} is also described in~\cite[\S 4.4]{Shao:2018:SPL}, where the authors show it is mathematically equivalent to Lanczos relationship~\eqref{eq:shao-lanczos}, with $\GruningU_k=\frac{1}{\sqrt{2}}U_k, \GruningW_k=\frac{1}{\sqrt{2}}V_k\GruningLowerBidiagonal_k^{-*}$, and $\GruningLowerBidiagonal_k$ the Cholesky factor of $T_k$. In particular, it is easy to see that relation~\eqref{eq:gruning-relationship} can be derived from~\eqref{eq:shao-lanczos}, using the similarity transformation
\[
\begin{bmatrix}
&T_k\\
I_k
\end{bmatrix}
=
\begin{bmatrix}
I_k&\\
   &\GruningLowerBidiagonal_k^{-*}
\end{bmatrix}
\begin{bmatrix}
                           &\GruningLowerBidiagonal_k\\
\GruningLowerBidiagonal_k^*&
\end{bmatrix}
\begin{bmatrix}
I_k&\\
   &\GruningLowerBidiagonal_k^*
\end{bmatrix},
\]
by substituting in~\eqref{eq:shao-lanczos} and multiplying on the right by 
$\frac{1}{\sqrt{2}}
\left[\begin{smallmatrix}
I_k&\\
   &\GruningLowerBidiagonal_k^{-*}
\end{smallmatrix}\right]
$.

In order to derive the algorithm, assume $\Gruningu_{j}, \GrHu_{j}, \Gruningw_{j-1}$ and $\Grbeta_{k+j-1}$ are available for a given $j\ge1$, with $\Gruningw_{0}=0$.
From column $j$ of \eqref{eq:gruning-relationship},
\[
\Grbeta_j\begin{bmatrix}\Gruningw_j\\-\bar{\Gruningw}_j\end{bmatrix}
=
H \begin{bmatrix}\Gruningu_j\\\bar{\Gruningu}_j\end{bmatrix}-
\Grbeta_{k+j-1}\begin{bmatrix}\Gruningw_{j-1}\\-\bar{\Gruningw}_{j-1}\end{bmatrix},
\]
thus,
\[
\Grbeta_j \Gruningw_j=\GrHu_j-\Grbeta_{k+j-1}\Gruningw_{j-1},
\]
whereby we can compute $\tilde{\Gruningw}_j=\Grbeta_j \Gruningw_j$, and then normalize to obtain $\Gruningw_j$ and $\GrHw_j$,
\[
\Gruningw_{j}=\tilde{\Gruningw}_{j}/\Grbeta_j,\;
\GrHw_j=x_j/\Grbeta_j,
\]
with $x_j=R\tilde{\Gruningw}_{j}-C\bar{\tilde{\Gruningw}}_{j}$ and $\Grbeta_j=\sqrt{2\text{Re}(\tilde{\Gruningw}_{j}^*x_{j})}$.

Similarly, from column $k+j$ of \eqref{eq:gruning-relationship},
\[
\Grbeta_{k+j} \Gruningu_{j+1}=\GrHw_j-\Grbeta_{j}\Gruningu_{j},
\]
which we can use to compute $\tilde{\Gruningu}_{j+1}=\Grbeta_{k+j} \Gruningu_{j+1}$, and normalize to obtain $\Gruningu_{j+1}$ and $\GrHu_{j+1}$,
\[
\Gruningu_{j+1}=\tilde{\Gruningu}_{j+1}/\Grbeta_{k+j},\;
\GrHu_{j+1}=y_{j+1}/\Grbeta_{k+j},
\]
with 
$y_{j+1}=R\tilde{\Gruningu}_{j+1}+C\bar{\tilde{\Gruningu}}_{j+1}$ and $\Grbeta_{k+j}=\sqrt{2\text{Re}(\tilde{\Gruningu}_{j+1}^*y_{j+1})}$.


An initial vector $\Gruningu_{1}$ can be generated from an arbitrary $\tilde{\Gruningu}_1$, by computing $y_{1}$ and normalizing,
\[
\Gruningu_{1}=\tilde{\Gruningu}_{1}/\Grbeta_0,
\qquad
\text{with }\Grbeta_0=\sqrt{2\text{Re}(\tilde{\Gruningu}_{1}^*y_{1})}.
\]

Further details of the method can be found in~\cite{Gruning:2011:ITL,Shao:2018:SPL}.

Bases $\GrHU_k$ and $\GrHW_k$ are constructed explicitly in our method as they are used in the reorthogonalization. Because of that, memory requirements for this method are larger than for the other two methods, as the number of vectors to keep doubles.

\subsubsection{Reorthogonalization}\label{sec:gruning-reorth}

Similarly to Shao's method, to accurately compute the eigenvectors it is important to guarantee the orthogonality of the generated vectors. This requires performing a reorthogonalization against the full basis after the local orthogonalization described above.

Assuming orthogonality of vector $\tilde{\Gruningu}_{j+1}$ is not perfect, we can write
\[
\Grbeta_{k+j}\begin{bmatrix}\Gruningu_{j+1}\\\bar{\Gruningu}_{j+1}\end{bmatrix}
=
\begin{bmatrix}\tilde{\Gruningu}_{j+1}\\\bar{\tilde{\Gruningu}}_{j+1}\end{bmatrix}
-\begin{bmatrix}\GruningU_j&\GruningW_{j}\\\mybar{\GruningU}_{j}&-\mybar{\GruningW}_{j}\end{bmatrix}
\begin{bmatrix}c\\d\end{bmatrix}.
\]
Premultiplying by
$\left[\begin{smallmatrix}
\GrHU_j       & \GrHW_j\\
-\bar{\GrHU}_j&\bar{\GrHW}_j
\end{smallmatrix}\right]^*\Hsign
$ and taking into account the orthogonality condition \eqref{eq:gruning-orthog},
\[
\begin{bmatrix}
c\\d
\end{bmatrix}
=
\begin{bmatrix}
\GrHU_j         & \GrHW_j\\
-\mybar{\GrHU}_j&\mybar{\GrHW}_j
\end{bmatrix}^*
\Hsign
\begin{bmatrix}\tilde{\Gruningu}_{j+1}\\\bar{\tilde{\Gruningu}}_{j+1}\end{bmatrix},
\]
which leads to
\[
\begin{aligned}
c&=2\text{Re}(\GrHU^*_j\tilde{\Gruningu}_{j+1}),\\
d&=2\text{Im}(\GrHW^*_j\tilde{\Gruningu}_{j+1})\text{i}.
\end{aligned}
\]
After computing $c,d$, a reorthogonalized vector $\check{\Gruningu}_{j+1}$ is obtained as
\[
\check{\Gruningu}_{j+1}=\Grbeta_{k+j}\Gruningu_{j+1}=\tilde{\Gruningu}_{j+1}- \GruningU_j c - \GruningW_j d.
\]

Our results indicate that ensuring orthogonality of each new vector $\Gruningu_{j+1}$ with respect to $\GruningU_j$ and $\GruningW_j$ is enough to also keep orthogonality of vectors $\Gruningw_j$. This is in line with what happens with the Lanczos bidiagonalization process \cite{Simon:2000:LMA}, and comes from the fact that vector bases $\GruningU_k$ and $\GruningW_k$ are closely related.

\subsubsection{Restart}
The restart must preserve the structure of the vector bases.
Let $\GruningLowerBidiagonal_k=\GruningQ\GruningLambda_k \GruningZ^*$ be the SVD of $\GruningLowerBidiagonal_k$, with $\GruningQ$ and $\GruningZ$ orthogonal. Then
\[
\begin{bmatrix}0&\GruningLowerBidiagonal_k\\ \GruningLowerBidiagonal_k^{*}&0\end{bmatrix}=
\begin{bmatrix}\GruningQ&0\\0&\GruningZ\end{bmatrix}
\begin{bmatrix}0&\GruningLambda_k\\\GruningLambda_k&0\end{bmatrix}
\begin{bmatrix}\GruningQ^*&0\\0&\GruningZ^*\end{bmatrix}.
\]
Substituting in \eqref{eq:gruning-relationship} and multiplying on the right by $\left[\begin{smallmatrix}\GruningQ&0\\0&\GruningZ\end{smallmatrix}\right]$ we get
\begin{equation}\label{eq:gruning-restart}
H
\begin{bmatrix}
\GruningU_k\GruningQ           &\GruningW_k\GruningZ\\
\mybar{\GruningU_k\GruningQ}&-\mybar{\GruningW_k\GruningZ}
\end{bmatrix}
=
\begin{bmatrix}
\GruningU_k\GruningQ           &\GruningW_k\GruningZ\\
\mybar{\GruningU_k\GruningQ}&-\mybar{\GruningW_k\GruningZ}
\end{bmatrix}
\begin{bmatrix}0&\GruningLambda_k\\\GruningLambda_k&0\end{bmatrix}
+\Grbeta_{2k}\begin{bmatrix}
\Gruningu_{k+1}\\\mybar{\Gruningu}_{k+1}
\end{bmatrix}
\begin{bmatrix}
0 & e_k^*\GruningZ
\end{bmatrix},
\end{equation}
or, defining $\hat{\GruningU}_k:=\GruningU_k Q, \hat{\GruningW}_k:=\GruningW_k P$ and $b := \Grbeta_{2k}\GruningZ^*e_k$, 
\begin{equation}
H
\begin{bmatrix}
\hat{\GruningU}_k           & \hat{\GruningW}_k\\
\mybar{\hat{\GruningU}}_k&-\mybar{\hat{\GruningW}}_k
\end{bmatrix}
=
\begin{bmatrix}
\hat{\GruningU}_k           &\hat{\GruningW}_k\\
\mybar{\hat{\GruningU}}_k&-\mybar{\hat{\GruningW}}_k
\end{bmatrix}
\begin{bmatrix}0&\GruningLambda_k\\\GruningLambda_k&0\end{bmatrix}
+\begin{bmatrix}
\Gruningu_{k+1}\\\mybar{\Gruningu}_{k+1}
\end{bmatrix}
\begin{bmatrix}
0 & b^*
\end{bmatrix}.
\label{eq:gruning-pre-truncated}
\end{equation}

\begin{algorithm}
  \caption{Gr\"uning's Lanczos with reorthogonalization}
  \label{alg:gruning-reorth-reset}
  \begin{algorithmic}[1]
    \REQUIRE A definite $2n\times 2n$ BSE matrix $H$ with structure \eqref{eq:bse1}; vector bases $\GruningU_{r+1}$, $\GrHU_{r+1}$, $\GruningW_r$ and $\GrHW_r$; bidiagonal $\GruningLowerBidiagonal_r$ with nonzero entries $\Grbeta_1,\dots,\Grbeta_r,\Grbeta_{k+1},\dots,\Grbeta_{k+r}$; vector $b$; number of Lanczos steps $k$.
    \ENSURE $\GruningU_{k+1},\GrHU_{k+1},\GruningW_k,\GrHW_k,\GruningLowerBidiagonal_k$, extending decomposition \eqref{eq:gruning-pre-truncated} to $2k$ columns.
    \FOR{$j=r+1,r+2,\dots,k$}
		\IF{$j=r+1$} 
			\STATE $\tilde{\Gruningw}_j = \GrHu_j - \GruningW_r b$
		\ELSE
			\STATE $\tilde{\Gruningw}_j = \GrHu_j - \Gruningw_{j-1}\Grbeta_{k+j-1}$
		\ENDIF
		\STATE $\tmpvec = R\tilde{\Gruningw}_j-C\bar{\tilde{\Gruningw}}_j$
		\STATE $\Grbeta_j = \sqrt{2\text{Re}(\tilde{\Gruningw}_{j}^*\tmpvec})$
		\STATE $\Gruningw_{j}=\tilde{\Gruningw}_{j}/\Grbeta_j,\quad \GrHw_{j}=\tmpvec/\Grbeta_j$
		\STATE $\tilde{\Gruningu}_{j+1} = \GrHw_{j} - \Gruningu_{j}\Grbeta_j$
		\STATE $c=2\text{Re}(\GrHU^*_j\tilde{\Gruningu}_{j+1}),\quad
				d=2\text{Im}(\GrHW^*_j\tilde{\Gruningu}_{j+1})\text{i},\quad
				\check{\Gruningu}_{j+1} = \tilde{\Gruningu}_{j+1} - \GruningU_j c - \GruningW_j d$
      	\STATE $\tmpvec = R\check{\Gruningu}_{j+1}+C\bar{\check{\Gruningu}}_{j+1}$
		\STATE $\Grbeta_{k+j} = \sqrt{2\text{Re}(\check{\Gruningu}_{j+1}^*\tmpvec})$
    	\STATE $\Gruningu_{j+1}=\check{\Gruningu}_{j+1}/\Grbeta_{k+j},\quad \GrHu_{j+1}=\tmpvec/\Grbeta_{k+j}$
    \ENDFOR
  \end{algorithmic}
\end{algorithm}

We assume the SVD decomposition of $\GruningLowerBidiagonal_k$ has been sorted so that the first elements in the diagonal of $\GruningLambda_k$ correspond to the smallest singular values. This decomposition is then truncated to $2r$ columns, with $r<k$, and then the process is restarted, extending the decomposition to $2k$ columns by means of Algorithm \ref{alg:gruning-reorth-reset}, which also details the reorthogonalization process described in \cref{sec:gruning-reorth}.
We repeat this process, extending and truncating until the desired number of Ritz eigenvalues have converged. 

The steps of the restart process are better visualized in~\cref{fig:gruning-restart}. The special orthogonalization with vector $b$ is considered in step 4. This can also be seen in~\cref{alg:gruning-reorth-reset}.

\def\mydiagonal#1{ 
  foreach \ii in {1,...,#1} { -- ++(0,-0.2) -- ++(0.2,0) }
  foreach \ii in {1,...,#1} { -- ++(0,0.2) -- ++(-0.2,0) }
}
\def\mylowerbidiagonalextra#1{ 
  -- ++(0,-0.2)
  foreach \ii in {1,...,#1} { -- ++(0,-0.2) -- ++(0.2,0) }
  -- ++(0,0.2)
  foreach \ii in {1,...,#1} { -- ++(0,0.2) -- ++(-0.2,0) }
}
\def\myupperbidiagonal#1#2{ 
  foreach \ii in {1,...,#1} { -- ++(0,-0.2) -- ++(0.2,0) }
  -- ++(0,0.2)
  foreach \ii in {1,...,#2} { -- ++(0,0.2) -- ++(-0.2,0) }
  -- ++(-0.2,0)
}

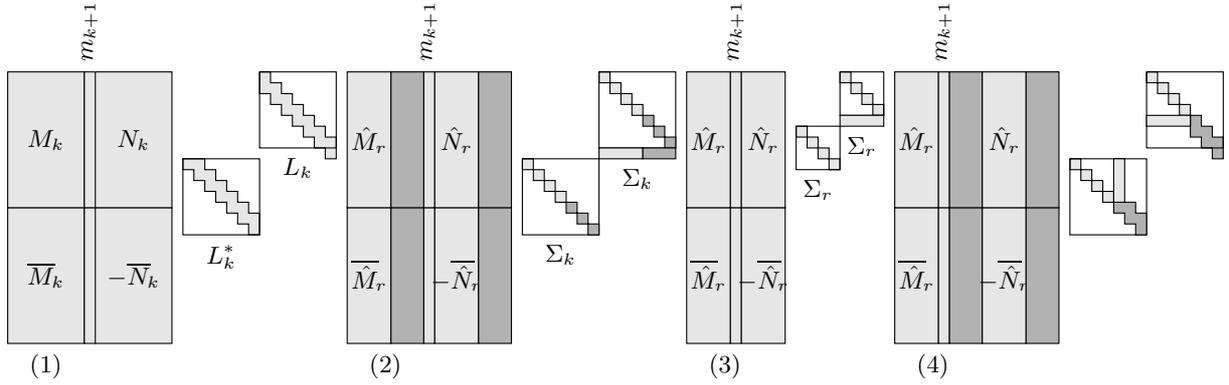
\begin{figure}
	\def\ligthgrey{black!10}
	\def\darkgrey{black!30}
  \resizebox{\textwidth}{!}{\begin{tikzpicture}[scale=0.7]
    \tikzstyle{every node}=[font=\footnotesize]
    \begin{scope}[xshift=0]
      \tikzmath{\nits=7;\startx=8.3;}
      \tikzmath{\startx=8.3;\movex=\nits*0.2;}
      \draw[fill=\ligthgrey] (\startx,0) rectangle node {$\mybar{\GruningU}_k$} +(\movex,5/2);
      \draw[fill=\ligthgrey] (\startx,5/2) rectangle node {$\GruningU_k$} +(\movex,5/2);
      \tikzmath{\startx=\startx+\movex;\movex=0.2;}
      \draw[fill=\ligthgrey] (\startx,0) rectangle +(\movex,5/2);
      \draw[fill=\ligthgrey] (\startx,5/2) rectangle +(\movex,5/2);
      \node[right,rotate=90] at (\startx+0.1,5) {$\Gruningu_{k+1}$};
      \tikzmath{\startx=\startx+\movex;\movex=\nits*0.2;}
      \draw[fill=\ligthgrey] (\startx,0) rectangle node {$-\mybar{\GruningW}_k$} +(\movex,5/2);
      \draw[fill=\ligthgrey] (\startx,5/2) rectangle node {$\GruningW_k$} +(\movex,5/2);
			\tikzmath{\startx=\startx+\movex+0.2;\starty=5-(\nits*0.2+0.2);}
			\tikzmath{\nitsmone=\nits-1;}
      \draw[fill=\ligthgrey] (\startx,\starty) \myupperbidiagonal{\nits}{\nitsmone};
			\tikzmath{\movex=\nits*0.2;}
      \node[below] at (\startx+\movex/2,\starty-\movex) {$\GruningLowerBidiagonal_k^*$};
      \draw (\startx,\starty) rectangle ++(\movex,-\movex);
			\tikzmath{\startx=\startx+\movex;\starty=5;\movex=0.2*\nits;}
      \draw[fill=\ligthgrey] (\startx,\starty) \mylowerbidiagonalextra{\nits};
      \node[below] at (\startx+\movex/2,\starty-\movex) {$\GruningLowerBidiagonal_k$};
			\tikzmath{\movex=\nits*0.2;}
      \draw (\startx,\starty) rectangle ++(\movex,-\movex);
			\tikzmath{\starty=5-(\nits*0.2+0.2);\movex=0.2*\nits;}
      \node[below,font=\footnotesize] at (9,0) {(1)};
    \end{scope}
    \begin{scope}[xshift=6.2cm]
      \tikzmath{\nits=4;\nitsb=3;\startx=8.3;}
      \tikzmath{\startx=8.3;\movex=\nits*0.2;}
      \draw[fill=\ligthgrey] (\startx,0) rectangle node {$\mybar{\hat{\GruningU}}_r$} +(\movex,5/2);
      \draw[fill=\ligthgrey] (\startx,5/2) rectangle node {$\hat{\GruningU}_r$} +(\movex,5/2);
      \tikzmath{\startx=\startx+\movex;\movex=\nitsb*0.2;}
      \draw[fill=\darkgrey] (\startx,0) rectangle +(\movex,5/2);
      \draw[fill=\darkgrey] (\startx,5/2) rectangle +(\movex,5/2);
      \tikzmath{\startx=\startx+\movex;\movex=0.2;}
      \draw[fill=\ligthgrey] (\startx,0) rectangle +(\movex,5/2);
      \draw[fill=\ligthgrey] (\startx,5/2) rectangle +(\movex,5/2);
      \node[right,rotate=90] at (\startx+0.1,5) {$\Gruningu_{k+1}$};
      \tikzmath{\startx=\startx+\movex;\movex=\nits*0.2;}
      \draw[fill=\ligthgrey] (\startx,0) rectangle node {$-\mybar{\hat{\GruningW}}_r$} +(\movex,5/2);
      \draw[fill=\ligthgrey] (\startx,5/2) rectangle node {$\hat{\GruningW}_r$} +(\movex,5/2);
      \tikzmath{\startx=\startx+\movex;\movex=\nitsb*0.2;}
      \draw[fill=\darkgrey] (\startx,0) rectangle +(\movex,5/2);
      \draw[fill=\darkgrey] (\startx,5/2) rectangle +(\movex,5/2);
			\tikzmath{\startx=\startx+\movex+0.2;\starty=5-(\nits*0.2+\nitsb*0.2+0.2);}
      \draw[fill=\ligthgrey] (\startx,\starty) \mydiagonal{\nits};
			\tikzmath{\movex=\nits*0.2+\nitsb*0.2;}
      \node[below] at (\startx+\movex/2,\starty-\movex) {$\GruningLambda_k$};
      \draw (\startx,\starty) rectangle ++(\movex,-\movex);
			\tikzmath{\movex=\nits*0.2;\startx=\startx+\movex;\starty=\starty-\movex;\movex=\nitsb*0.2;}
      \draw[fill=\darkgrey] (\startx,\starty) \mydiagonal{\nitsb};
			\tikzmath{\startx=\startx+\movex;\starty=5;\movex=0.2*\nits;}
      \draw[fill=\ligthgrey] (\startx,\starty) \mydiagonal{\nits};
			\tikzmath{\movex=\nits*0.2+\nitsb*0.2;}
      \draw (\startx,\starty) rectangle ++(\movex,-\movex);
      \draw[fill=\darkgrey] (\startx+\nits*0.2,\starty-\nits*0.2) \mydiagonal{\nitsb};
      \node[below] at (\startx+\movex/2,\starty-\movex-0.2) {$\GruningLambda_k$};
			\tikzmath{\starty=5-(\nits*0.2+\nitsb*0.2+0.2);\movex=0.2*\nits;}
      \draw[fill=\ligthgrey] (\startx,\starty) rectangle +(\movex,0.2);
			\tikzmath{\startx=\startx+\movex;\starty=\starty;\movex=0.2*\nitsb;}
      \draw[fill=\darkgrey] (\startx,\starty) rectangle +(\movex,0.2);
      \node[below,font=\footnotesize] at (9,0) {(2)};
    \end{scope}
    \begin{scope}[xshift=12.4cm]
      \tikzmath{\nits=4;\startx=8.3;}
      \tikzmath{\startx=8.3;\movex=\nits*0.2;}
      \draw[fill=\ligthgrey] (\startx,0) rectangle node {$\mybar{\hat{\GruningU}}_r$} +(\movex,5/2);
      \draw[fill=\ligthgrey] (\startx,5/2) rectangle node {$\hat{\GruningU}_r$} +(\movex,5/2);
      \tikzmath{\startx=\startx+\movex;\movex=0.2;}
      \draw[fill=\ligthgrey] (\startx,0) rectangle +(\movex,5/2);
      \draw[fill=\ligthgrey] (\startx,5/2) rectangle +(\movex,5/2);
      \node[right,rotate=90] at (\startx+0.1,5) {$\Gruningu_{k+1}$};
      \tikzmath{\startx=\startx+\movex;\movex=\nits*0.2;}
      \draw[fill=\ligthgrey] (\startx,0) rectangle node {$-\mybar{\hat{\GruningW}}_r$} +(\movex,5/2);
      \draw[fill=\ligthgrey] (\startx,5/2) rectangle node {$\hat{\GruningW}_r$} +(\movex,5/2);
			\tikzmath{\startx=\startx+\movex+0.2;\starty=5-(\nits*0.2+0.2);}
      \draw[fill=\ligthgrey] (\startx,\starty) \mydiagonal{\nits};
      \node[below] at (\startx+\movex/2,\starty-\movex) {$\GruningLambda_r$};
			\tikzmath{\movex=\nits*0.2;}
      \draw (\startx,\starty) rectangle ++(\movex,-\movex);
			\tikzmath{\startx=\startx+\movex;\starty=5;\movex=0.2*\nits;}
      \draw[fill=\ligthgrey] (\startx,\starty) \mydiagonal{\nits};
      \node[below] at (\startx+\movex/2,\starty-\movex-0.2) {$\GruningLambda_r$};
			\tikzmath{\movex=\nits*0.2;}
      \draw (\startx,\starty) rectangle ++(\movex,-\movex);
			\tikzmath{\starty=5-(\nits*0.2+0.2);\movex=0.2*\nits;}
      \draw[fill=\ligthgrey] (\startx,\starty) rectangle +(\movex,0.2);
      \node[below,font=\footnotesize] at (9,0) {(3)};
    \end{scope}
    \begin{scope}[xshift=16.2cm]
      \tikzmath{\nits=4;\nitsb=3;\startx=8.3;}
      \tikzmath{\startx=8.3;\movex=\nits*0.2;}
      \draw[fill=\ligthgrey] (\startx,0) rectangle node {$\mybar{\hat{\GruningU}}_r$} +(\movex,5/2);
      \draw[fill=\ligthgrey] (\startx,5/2) rectangle node {$\hat{\GruningU}_r$} +(\movex,5/2);
      \tikzmath{\startx=\startx+\movex;\movex=0.2;}
      \draw[fill=\ligthgrey] (\startx,0) rectangle +(\movex,5/2);
      \draw[fill=\ligthgrey] (\startx,5/2) rectangle +(\movex,5/2);
      \node[right,rotate=90] at (\startx+0.1,5) {$\Gruningu_{k+1}$};
      \tikzmath{\startx=\startx+\movex;\movex=\nitsb*0.2;}
      \draw[fill=\darkgrey] (\startx,0) rectangle +(\movex,5/2);
      \draw[fill=\darkgrey] (\startx,5/2) rectangle +(\movex,5/2);
      \tikzmath{\startx=\startx+\movex;\movex=\nits*0.2;}
      \draw[fill=\ligthgrey] (\startx,0) rectangle node {$-\mybar{\hat{\GruningW}}_r$} +(\movex,5/2);
      \draw[fill=\ligthgrey] (\startx,5/2) rectangle node {$\hat{\GruningW}_r$} +(\movex,5/2);
      \tikzmath{\startx=\startx+\movex;\movex=\nitsb*0.2;}
      \draw[fill=\darkgrey] (\startx,0) rectangle +(\movex,5/2);
      \draw[fill=\darkgrey] (\startx,5/2) rectangle +(\movex,5/2);
			\tikzmath{\startx=\startx+\movex+0.2;\starty=5-(\nits*0.2+\nitsb*0.2+0.2);}
      \draw[fill=\ligthgrey] (\startx,\starty) \mydiagonal{\nits};
			\tikzmath{\movex=\nits*0.2+\nitsb*0.2;}
      \draw (\startx,\starty) rectangle ++(\movex,-\movex);
			\tikzmath{\movex=\nits*0.2;\startx=\startx+\movex;\starty=\starty-\movex;\movex=\nitsb*0.2;}
			\tikzmath{\nitsbmone=\nitsb-1;}
      \draw[fill=\darkgrey] (\startx,\starty) \myupperbidiagonal{\nitsb}{\nitsbmone};
      \draw[fill=\ligthgrey] (\startx,\starty) rectangle +(0.2,\nits*0.2);
			\tikzmath{\startx=\startx+\movex;\starty=5;\movex=0.2*\nits;}
      \draw[fill=\ligthgrey] (\startx,\starty) \mydiagonal{\nits};
      \draw[fill=\darkgrey] (\startx+\nits*0.2,\starty-\nits*0.2) \mylowerbidiagonalextra{\nitsb};
			\tikzmath{\movex=\nits*0.2+\nitsb*0.2;}
      \draw (\startx,\starty) rectangle ++(\movex,-\movex);
			\tikzmath{\starty=5-(\nits*0.2+0.2);\movex=0.2*\nits;}
      \draw[fill=\ligthgrey] (\startx,\starty) rectangle +(\movex,0.2);
      \node[below,font=\footnotesize] at (9,0) {(4)};
    \end{scope}
  \end{tikzpicture}}
  \caption{Illustration of the  steps of thick restart in the Gr\"uning method: (1) initial Lanczos factorization of order $k$, (2) solve projected problem, sort and check convergence, (3) truncate to factorization of order $r$, and (4) extend to a factorization of order $k$.}
  \label{fig:gruning-restart}
\end{figure}

\subsubsection{Obtaining the eigenvectors}
To obtain the eigenvalues and eigenvectors, we take the Krylov factorization~\eqref{eq:gruning-pre-truncated} and define 
\begin{equation}\label{eq:GruningP}
\GruningP:=\frac{1}{\sqrt{2}}\begin{bmatrix}I_k&I_k\\I_k&-I_k\end{bmatrix},
\end{equation}
so that
\[
\begin{bmatrix}0&\GruningLambda_k\\\GruningLambda_k&0\end{bmatrix}=
\GruningP\begin{bmatrix}\GruningLambda_k&0\\0&-\GruningLambda_k\end{bmatrix}\GruningP^*.
\]
Substituting this equality in~\eqref{eq:gruning-pre-truncated} and multiplying on the right by $\GruningP$ yields
\[
H
\begin{bmatrix}
\tilde{X}_1 & \mybar{\tilde{X}}_2 \\
\tilde{X}_2 & \mybar{\tilde{X}}_1 \\
\end{bmatrix}
=
\begin{bmatrix}
\tilde{X}_1 & \mybar{\tilde{X}}_2 \\
\tilde{X}_2 & \mybar{\tilde{X}}_1 \\
\end{bmatrix}
\begin{bmatrix}
\tilde{\Lambda}_+& 0\\
0                &-\tilde{\Lambda}_+
\end{bmatrix}
+
\frac{1}{\sqrt{2}}
\begin{bmatrix}
\Gruningu_{k+1}\\\mybar{\Gruningu}_{k+1}
\end{bmatrix}
\begin{bmatrix}
b^* & -b^*
\end{bmatrix},
\]
where we use the following definitions, corresponding to the approximate eigenpairs:
\[
\begin{bmatrix}
\tilde{X}_1 & \mybar{\tilde{X}}_2 \\
\tilde{X}_2 & \mybar{\tilde{X}}_1 \\
\end{bmatrix}
:=
\begin{bmatrix}
\hat{\GruningU}_k           &\hat{\GruningW}_k\\
\mybar{\hat{\GruningU}}_k&-\mybar{\hat{\GruningW}}_k
\end{bmatrix}
\GruningP
=
\frac{1}{\sqrt{2}}
\begin{bmatrix}
\hat{\GruningU}_k + \hat{\GruningW}_k  &
\hat{\GruningU}_k - \hat{\GruningW}_k  \\
\mybar{\hat{\GruningU}}_k - \mybar{\hat{\GruningW}}_k & 
\mybar{\hat{\GruningU}}_k + \mybar{\hat{\GruningW}}_k
\end{bmatrix},
\quad
\tilde{\Lambda}_+:=\GruningLambda_k.
\]
The product with $\frac{1}{\sqrt{2}}$ is not required, because the eigenvectors can be normalized in different ways. The normalization of the vectors will be discussed in~\cref{sec:implem}. The left eigenvectors are computed as in~\eqref{eq:spliteigenvectors}.

\subsection{Method resulting in a Bethe--Salpeter projected eigenproblem}\label{sec:projectedbse}
\def\ProjBSEX{W}
\def\ProjBSEY{Z}
\def\ProjBSEx{w}
\def\ProjBSEy{z}
\def\ProjBSEA{A}
\def\ProjBSEB{B}
\def\ProjBSEa{a}
\def\ProjBSEb{b}
\def\ProjBSEZ{\tilde{G}}

\subsubsection{Lanczos recurrence}

In~\cite[\S 4.4]{Shao:2018:SPL}, the authors mention a reformulation of the Lanczos relation where the projected matrix also has Bethe--Salpeter structure. Here we work out the details of a restarted solver based on this formulation. We start by stating a proposition showing the equivalence between the two Lanczos decompositions~\cite[\S 4.4]{Shao:2018:SPL}.

\begin{proposition}\label{prop:projectedbse}
The Lanczos-type relation
\begin{equation}
H
\begin{bmatrix}
\ProjBSEX_k & \mybar{\ProjBSEY}_k\\
\ProjBSEY_k & \mybar{\ProjBSEX}_k
\end{bmatrix}
=
\begin{bmatrix}
\ProjBSEX_k & \mybar{\ProjBSEY}_k\\
\ProjBSEY_k & \mybar{\ProjBSEX}_k
\end{bmatrix}
\begin{bmatrix}
\ProjBSEA_k & \ProjBSEB_k\\
-\ProjBSEB_k & -\ProjBSEA_k
\end{bmatrix}
+
\frac{1}{2}\beta_k
\begin{bmatrix}
\ProjBSEx_{k+1}&\bar{\ProjBSEy}_{k+1}\\
\ProjBSEy_{k+1}&\bar{\ProjBSEx}_{k+1}
\end{bmatrix}
\begin{bmatrix}
e_{k}^*&-e_{k}^*\\
e_{k}^*&-e_{k}^*
\end{bmatrix},
\label{eq:projectedbse}
\end{equation}
where
$$
\ProjBSEX_{k+1}:=\frac{U_{k+1}+V_{k+1}}{2},\quad \ProjBSEY_{k+1}:=\frac{\mybar{U}_{k+1}-\mybar{V}_{k+1}}{2},\quad
\ProjBSEA_k:=\frac{I+T_k}{2},\quad \ProjBSEB_k:=\frac{I-T_k}{2},
$$
is equivalent to~\eqref{eq:shao-lanczos} in the sense that both provide the same Ritz approximations.
\end{proposition}
\begin{proof}
The proof is based on the similarity transformation
\[
\begin{bmatrix}
&T_k\\
I_k
\end{bmatrix}
=
\GruningP
\begin{bmatrix}
A_k&B_k\\
-B_k&-A_k
\end{bmatrix}
\GruningP,
\]
where 
$\GruningP=\GruningP^{-1}=\frac{1}{\sqrt{2}}
\left[\begin{smallmatrix}
I_k&I_k\\I_k&-I_k
\end{smallmatrix}\right]
$ as defined in~\eqref{eq:GruningP}.

Substituting in~\eqref{eq:shao-lanczos} and multiplying on the right by $\frac{1}{\sqrt{2}}\GruningP$, we get the relation \eqref{eq:projectedbse}\footnote{The last term in~\eqref{eq:projectedbse} differs from the one given in~\cite[\S 4.4]{Shao:2018:SPL}, which we believe is wrong.}. It follows that the Ritz approximations provided by~\eqref{eq:projectedbse} and~\eqref{eq:shao-lanczos} are the same.
\end{proof}

From the relation with the vector bases of Shao's method, we can easily check that the orthogonality condition is
\begin{equation}
\begin{bmatrix}
\ProjBSEX_k  &-\mybar{\ProjBSEY}_k\\
-\ProjBSEY_k & \mybar{\ProjBSEX}_k
\end{bmatrix}^*
\begin{bmatrix}
\ProjBSEX_k & \mybar{\ProjBSEY}_k\\
\ProjBSEY_k & \mybar{\ProjBSEX}_k
\end{bmatrix}
=
I_{2k},
\label{eq:projectedbse-ortog}
\end{equation}
or
\begin{equation}
\begin{bmatrix}
\ProjBSEX_k  &\mybar{\ProjBSEY}_k\\
\ProjBSEY_k & \mybar{\ProjBSEX}_k
\end{bmatrix}^*
\begin{bmatrix}
I_n&\\&-I_n
\end{bmatrix}
\begin{bmatrix}
\ProjBSEX_k & \mybar{\ProjBSEY}_k\\
\ProjBSEY_k & \mybar{\ProjBSEX}_k
\end{bmatrix}
=
\begin{bmatrix}
I_k&\\&-I_k
\end{bmatrix}.
\end{equation}

To derive the algorithm, we have to consider~\eqref{eq:projectedbse} as a kind of block Lanczos-type relation with two initial vectors. We use the notation
\begin{align*}
\ProjBSEA_k&=\text{tridiag}
\left\lbrace
\begingroup 
\setlength\arraycolsep{2pt}
\begin{matrix}
        &\beta_1/2&      &\cdots&      &\beta_{k-1}/2& \\
\ProjBSEa_1&       &\cdots&      &\cdots&           & \ProjBSEa_{k}\\
        &\beta_1/2&      &\cdots&      &\beta_{k-1}/2&
\end{matrix}
\endgroup
\right\rbrace,\\
\ProjBSEB_k&=\text{tridiag}
\left\lbrace
\begingroup 
\setlength\arraycolsep{2pt}
\begin{matrix}
        &-\beta_1/2&      &\cdots&      &-\beta_{k-1}/2& \\
\ProjBSEb_1&       &\cdots&      &\cdots&           & \ProjBSEb_{k}\\
        &-\beta_1/2&      &\cdots&      &-\beta_{k-1}/2&
\end{matrix}
\endgroup
\right\rbrace,
\end{align*}
with $\ProjBSEa_i=(1+\alpha_i)/2$, $\ProjBSEb_i=(1-\alpha_i)/2$, $i=1,\dots,k$.

In each iteration $j$ of the method we start with vectors $u_j, v_j$ corresponding to Shao's method, as well as vectors $\ProjBSEx_{j-1},\ProjBSEy_{j-1},\ProjBSEx_{j}$ and $\ProjBSEy_{j}$, and we compute vectors $u_{j+1}, v_{j+1}$ and $\ProjBSEx_{j+1},\ProjBSEy_{j+1}$. Vectors $u_j,v_j$ are discarded once the iteration is complete.

Equating column $j$ ($j=2,\dots,k$) in both sides of~\eqref{eq:projectedbse} results in
\begin{equation}
\begin{bmatrix}\ProjBSEx'_j\\\ProjBSEy'_j\end{bmatrix}
=
\frac{\beta_{j-1}}{2}
\begin{bmatrix}
\ProjBSEx_{j-1}+\overline{\ProjBSEy}_{j-1}\\
\ProjBSEy_{j-1}+\overline{\ProjBSEx}_{j-1}
\end{bmatrix}
+\ProjBSEa_j
\begin{bmatrix}\ProjBSEx_j\\\ProjBSEy_j\end{bmatrix}
+(\ProjBSEa_j-1)
\begin{bmatrix}\overline{\ProjBSEy}_j\\\overline{\ProjBSEx}_j\end{bmatrix}
+
\frac{\beta_j}{2}
\begin{bmatrix}
\ProjBSEx_{j+1}+\overline{\ProjBSEy}_{j+1}\\
\overline{\ProjBSEx}_{j+1}+\ProjBSEy_{j+1}
\end{bmatrix},
\label{eq:projectedbse-col}
\end{equation}
where
\[
\begin{bmatrix}\ProjBSEx'_j\\\ProjBSEy'_j\end{bmatrix}=
H\begin{bmatrix}\ProjBSEx_j\\\ProjBSEy_j\end{bmatrix}.
\]
Note that equating column $j+k$ gives the same vectors, but swapped, conjugated and with the sign changed.
Multiplying on the left by $[\ProjBSEx_j^*, -\ProjBSEy_j^*]$ and taking into account orthogonality condition \eqref{eq:projectedbse-ortog}, yields
\[
\ProjBSEa_j=\begin{bmatrix}\ProjBSEx_j\\-\ProjBSEy_j\end{bmatrix}^*\begin{bmatrix}\ProjBSEx'_j\\\ProjBSEy'_j\end{bmatrix}.
\]
Vectors $\ProjBSEx'_j$ and $\ProjBSEy'_j$ can be computed efficiently using $v_j$. Equation \eqref{eq:shao-lanczos} from Shao's method implies
\begin{equation}
\begin{bmatrix}v_j\\-\bar{v}_j\end{bmatrix}
=
H\begin{bmatrix}u_j\\\bar{u}_j\end{bmatrix}.
\label{eq:shao-relation-u-v}
\end{equation}
Using that, we get
\[
\begin{bmatrix}\ProjBSEx'_j\\\ProjBSEy'_j\end{bmatrix}
=
\frac{1}{2}H\begin{bmatrix}u_j+v_j\\\overline{u}_j-\overline{v}_j\end{bmatrix}
=
\frac{1}{2}\left(
\begin{bmatrix}v_j\\-\overline{v}_j\end{bmatrix}+
H\begin{bmatrix}v_j\\-\overline{v}_j\end{bmatrix}
\right).
\]
Thus,
\[
\ProjBSEx'_j=\frac{1}{2}(v_j+v_j'),\quad \ProjBSEy'_j=\frac{1}{2}(\overline{v_j'-v_j}), 
\]
where
\[
v_j'=Rv_j-C\overline{v}_j.
\]

In the relation \eqref{eq:projectedbse-col} there are two unknowns, $\ProjBSEx_{j+1}$ and $\ProjBSEy_{j+1}$. If we express them in terms of $u_{k+1}$, it is possible to isolate them,
\[
\frac{1}{2}\beta_j
\begin{bmatrix}
u_{j+1}\\
\overline{u}_{j+1}
\end{bmatrix}
=
\begin{bmatrix}
\ProjBSEx'_j
-\frac{1}{2}\beta_{j-1}(\ProjBSEx_{j-1}+\overline{\ProjBSEy}_{j-1})-\ProjBSEa_j\ProjBSEx_j
-(\ProjBSEa_j-1)\overline{\ProjBSEy}_j
\\
\ProjBSEy'_j
-\frac{1}{2}\beta_{j-1}(\ProjBSEy_{j-1}+\overline{\ProjBSEx}_{j-1})-\ProjBSEa_j\ProjBSEy_j
-(\ProjBSEa_j-1)\overline{\ProjBSEx}_j
\end{bmatrix}.
\]

We only need the upper block of the equality to compute $\tilde{u}_{j+1}=\frac{1}{2}\beta_j u_{j+1}$.
From that, we compute $\beta_j, u_{j+1}, v_{j+1}$, using \eqref{eq:shao-relation-u-v} and taking into account that, from the formulation of Shao's method, we have
\[
\begin{bmatrix}
u_{j+1}\\-\bar{u}_{j+1}
\end{bmatrix}^*
\begin{bmatrix}
v_{j+1}\\-\bar{v}_{j+1}
\end{bmatrix}
=2.
\]
If we compute
\[
\tilde{v}_{j+1}
=R\tilde{u}_{j+1}+C\overline{\tilde{u}}_{j+1},
\]
we have that
\[
\begin{bmatrix}
\tilde{u}_{j+1}\\-\tilde{\bar{u}}_{j+1}
\end{bmatrix}^*
\begin{bmatrix}
\tilde{v}_{j+1}\\-\tilde{\bar{v}}_{j+1}
\end{bmatrix}
=\frac{1}{2}\beta_j^2
\implies
\beta_j=2\sqrt{\text{Re}(\tilde{u}_{j+1}^*\tilde{v}_{j+1})},
\]
$$
u_{j+1}=\frac{2}{\beta_j}\tilde{u}_{j+1},
\quad
v_{j+1}=\frac{2}{\beta_j}\tilde{v}_{j+1},
$$
and finally we compute the next vectors $\ProjBSEx_{j+1}$, $\ProjBSEy_{j+1}$, which ends the iteration,
$$
\ProjBSEx_{j+1}=\frac{u_{j+1}+v_{j+1}}{2},\quad \ProjBSEy_{j+1}=\frac{\bar{u}_{j+1}-\bar{v}_{j+1}}{2}.
$$

\subsubsection{Reorthogonalization}
To derive the reorthogonalization coefficients, we consider the equality,
\[
\frac{1}{2}\beta_j
\begin{bmatrix}
u_{j+1}\\
\overline{u}_{j+1}
\end{bmatrix}
=
\begin{bmatrix}
\tilde{u}_{j+1}\\
\overline{\tilde{u}}_{j+1}
\end{bmatrix}
-
\begin{bmatrix}
\ProjBSEX_j & \mybar{\ProjBSEY}_j\\
\ProjBSEY_j & \mybar{\ProjBSEX}_j
\end{bmatrix}
\begin{bmatrix}
c\\d
\end{bmatrix}.
\]
Premultiplying by
$$
\begin{bmatrix}
\ProjBSEX_j  &-\mybar{\ProjBSEY}_j\\
-\ProjBSEY_j & \mybar{\ProjBSEX}_j
\end{bmatrix}^*
$$
and taking into account the orthogonality condition~\eqref{eq:projectedbse-ortog}, we have
\begin{equation}
\begin{bmatrix}
c\\
d
\end{bmatrix}
=
\begin{bmatrix}
\ProjBSEX_j  &-\mybar{\ProjBSEY}_j\\
-\ProjBSEY_j & \mybar{\ProjBSEX}_j
\end{bmatrix}^*
\begin{bmatrix}
\tilde{u}_{j+1}\\
\overline{\tilde{u}}_{j+1}
\end{bmatrix},
\label{eq-ortog-coefs}
\end{equation}
resulting in,
\[
c=\bar{d}=\ProjBSEX_j \tilde{u}_{j+1} - \mybar{\ProjBSEY}_j \overline{\tilde{u}}_{j+1}.
\]
We can then obtain a reorthogonalized vector $\check{u}_{j+1}$ as
\[
\check{u}_{j+1}=\frac{1}{2}\beta_j u_{j+1} = \tilde{u}_{j+1}- \ProjBSEX_j c - \overline{\ProjBSEY}_j \bar{c}.
\]

\begin{algorithm}
  \caption{Projected BSE method with reorthogonalization}
  \label{alg:projbse-reorth}
  \begin{algorithmic}[1]
    \REQUIRE A definite $2n\times 2n$ BSE matrix $H$ with structure \eqref{eq:bse1}; initial vectors in columns of $\ProjBSEX_{r+1}$ and $\ProjBSEY_{r+1}$; vector $b$; number of Lanczos steps $k$.
    \ENSURE $\ProjBSEX_k,\ProjBSEY_k,T_k$ making up Lanczos-type decomposition \eqref{eq:projectedbse}.
    \STATE $v=\ProjBSEx_{r+1}-\overline{\ProjBSEy}_{r+1}$
    \FOR{$j=r+1,r+2,\dots,k$}
    \STATE $v'=Rv-C\bar{v}$
    \STATE $\ProjBSEx'=(v'+v)/2,\quad \ProjBSEy'=(\overline{v'-v})/2$
    \STATE $\tilde{\ProjBSEa}_j=\ProjBSEx_j^*\ProjBSEx'-\ProjBSEy_j^*\ProjBSEy'$
    \IF {$j=r+1$}
   	\STATE $\tilde{u}=\ProjBSEx' - \ProjBSEX_r b-\overline{\ProjBSEY}_r b - 
   			\tilde{\ProjBSEa}_j\ProjBSEx_j - (\tilde{\ProjBSEa}_j-1)\overline{\ProjBSEy}_j$
    \ELSE
	\STATE $\tilde{u}=\ProjBSEx'-\frac{1}{2}\beta_{j-1}(\ProjBSEx_{j-1}+\overline{\ProjBSEy}_{j-1})
    		-\tilde{\ProjBSEa}_j\ProjBSEx_j-(\tilde{\ProjBSEa}_j-1)\overline{\ProjBSEy}_j$
    \ENDIF
    \STATE $c=\ProjBSEX_j^*\tilde{u}-\ProjBSEY_j^*\overline{\tilde{u}}$
    \STATE $\check{u}=\tilde{u}-\ProjBSEX_j c-\overline{\ProjBSEY_j c}$
    \STATE $\ProjBSEa_j = \tilde{\ProjBSEa}_j+c_j,\quad \alpha_j=2a_j-1$
    \STATE $\tilde{v}=R\check{u}+C\bar{\check{u}}$
    \STATE $\beta_j=2\sqrt{\text{Re}(\check{u}^*\tilde{v})}$
    \STATE $u=2/\beta_j\check{u},\quad v=2/\beta_j\tilde{v}$
    \STATE $\ProjBSEx_{j+1}=(u+v)/2,\quad \ProjBSEy_{j+1}=(\overline{u-v})/2$
    \ENDFOR
  \end{algorithmic}
\end{algorithm}

\subsubsection{Restart}

The projected matrix $\left[\begin{smallmatrix}\ProjBSEA_k & \ProjBSEB_k\\ -\ProjBSEB_k & -\ProjBSEA_k\end{smallmatrix}\right]$ is real and has a Bethe--Salpeter structure. To restart the Lanczos relation~\eqref{eq:projectedbse}, we can exploit the fact that both $\ProjBSEA_k$ and $\ProjBSEB_k$ depend on $T_k$, and use its eigendecomposition  $T_k=QD_kQ^*$, with $Q$ orthogonal. Then,
$$
\begin{bmatrix}
\ProjBSEA_k & \ProjBSEB_k\\
-\ProjBSEB_k & -\ProjBSEA_k
\end{bmatrix}
=
\frac{1}{2}
\begin{bmatrix}
 I+T_k & I-T_k\\
-I+T_k &-I-T_k
\end{bmatrix}
=
\begin{bmatrix}
Q & \\
  &Q
\end{bmatrix}
\frac{1}{2}
\begin{bmatrix}
 I+D_k  & I-D_k\\
-(I-D_k)&-(I+D_k)
\end{bmatrix}
\begin{bmatrix}
Q^* & \\
    &Q^*
\end{bmatrix}.
$$
Inserting this into~\eqref{eq:projectedbse} and multiplying on the right by $\left[\begin{smallmatrix} Q &0 \\ 0&Q \end{smallmatrix}\right]$ we arrive at
\begin{equation}\label{eq:restart-projbse}
H
\begin{bmatrix}
\ProjBSEX_kQ & \mybar{\ProjBSEY}_kQ\\
\ProjBSEY_kQ & \mybar{\ProjBSEX}_kQ
\end{bmatrix}
=
\begin{bmatrix}
\ProjBSEX_kQ & \mybar{\ProjBSEY}_kQ\\
\ProjBSEY_kQ & \mybar{\ProjBSEX}_kQ
\end{bmatrix}
\frac{1}{2}
\begin{bmatrix}
 I+D_k  & I-D_k\\
-(I-D_k)&-(I+D_k)
\end{bmatrix}
+
\frac{1}{2}\beta_k
\begin{bmatrix}
u_{k+1}\\\bar{u}_{k+1}
\end{bmatrix}
\begin{bmatrix}
e_{k}^*Q&-e_{k}^*Q
\end{bmatrix},
\end{equation}
or, defining $\tilde{\ProjBSEX}_k:=\ProjBSEX_kQ$, $\tilde{\ProjBSEY}_k:=\ProjBSEY_kQ$,
\begin{equation}
H
\begin{bmatrix}
\tilde{\ProjBSEX}_k & \mybar{\tilde{\ProjBSEY}}_k\\
\tilde{\ProjBSEY}_k & \mybar{\tilde{\ProjBSEX}}_k
\end{bmatrix}
=
\begin{bmatrix}
\tilde{\ProjBSEX}_k & \mybar{\tilde{\ProjBSEY}}_k\\
\tilde{\ProjBSEY}_k & \mybar{\tilde{\ProjBSEX}}_k
\end{bmatrix}
\frac{1}{2}
\begin{bmatrix}
 I+D_k  & I-D_k\\
-(I-D_k)&-(I+D_k)
\end{bmatrix}
+
\frac{1}{2}\beta_k
\begin{bmatrix}
u_{k+1}\\\bar{u}_{k+1}
\end{bmatrix}
\begin{bmatrix}
e_{k}^*Q&-e_{k}^*Q
\end{bmatrix}.
\label{eq:projectedbse-restart}
\end{equation}
Similarly to the other two methods, this expression can be truncated to $2r$ columns, after reordering to keep the wanted approximate eigenvalues. Then, it can be expanded again to $2k$ columns by means of Algorithm \ref{alg:projbse-reorth}.

\subsubsection{Obtaining the eigenvectors}

It can be easily verified that the eigenvalue decomposition of the projected matrix can be written as
$$
\frac{1}{2}
\begin{bmatrix}
 I+D_k  & I-D_k\\
-(I-D_k)&-(I+D_k)
\end{bmatrix}
\ProjBSEZ
=
\ProjBSEZ
\begin{bmatrix}
D_k^\frac{1}{2}&0\\
0       &-D_k^\frac{1}{2}
\end{bmatrix},
$$
where
\begin{equation*}
\ProjBSEZ:=\begin{bmatrix}
D_k^\frac{1}{2}+I&D_k^\frac{1}{2}-I\\
D_k^\frac{1}{2}-I&D_k^\frac{1}{2}+I
\end{bmatrix}.
\end{equation*}
With this, if we multiply the relation~\eqref{eq:projectedbse-restart} on the right by $\ProjBSEZ$, we obtain
\begin{equation}
H
\begin{bmatrix}
\tilde{X}_1 & \mybar{\tilde{X}}_2 \\
\tilde{X}_2 & \mybar{\tilde{X}}_1 \\
\end{bmatrix}
=
\begin{bmatrix}
\tilde{X}_1 & \mybar{\tilde{X}}_2 \\
\tilde{X}_2 & \mybar{\tilde{X}}_1 \\
\end{bmatrix}
\begin{bmatrix}
\tilde\Lambda_+  & \\
          &-\tilde\Lambda_+
\end{bmatrix}
+
\frac{1}{2}\beta_k
\begin{bmatrix}
u_{k+1}\\\bar{u}_{k+1}
\end{bmatrix}
\begin{bmatrix}
e_{k}^*Q\ProjBSEZ&-e_{k}^*Q\ProjBSEZ
\end{bmatrix},
\end{equation}
where we have defined the approximate eigenpairs as
$$
\begin{bmatrix}
\tilde{X}_1 & \mybar{\tilde{X}}_2 \\
\tilde{X}_2 & \mybar{\tilde{X}}_1 \\
\end{bmatrix}
:=
\begin{bmatrix}
\tilde{\ProjBSEX}_k & \mybar{\tilde{\ProjBSEY}}_k\\
\tilde{\ProjBSEY}_k & \mybar{\tilde{\ProjBSEX}}_k
\end{bmatrix}
\ProjBSEZ,
\qquad
\tilde{\Lambda}_+:=D_k^\frac{1}{2}.
$$

\section{A more general view of the decompositions used}\label{sec:unified}
In this section, we analyze the three Lanczos decompositions used in the previously presented methods, providing a more general form for each of them, and presenting forms for similiarity transformations that preserve the structure of each decomposition.

In order to reach a more general description of the decompositions, we can resort to the concept of Krylov decompositions introduced by Stewart in his seminal paper~\cite{Stewart:2001:KAL}. Recall that, given a square matrix $A$, a Krylov decomposition of order $k$ is a relation of the form $AU_k=U_kC_k+u_{k+1}c_{k+1}^*$, where $C_k$ is of order $k$ and the columns of $[U_k,u_{k+1}]$ are linearly independent. It is a generalization of the Arnoldi decomposition, where the basis of the decomposition (the columns of $[U_k,u_{k+1}]$, which span a Krylov subspace for $A$) are not necessarily orthonormal, the Rayleigh quotient $C_k$ is not restricted to an upper Hessenberg form, and the rank-1 residual term $u_{k+1}c_{k+1}^*$ is not limited to the $k$th column. Stewart~\cite{Stewart:2001:KAL} uses two types of transformations, similarity and translation, that preserve the form of a Krylov decomposition, and can be used to manipulate it in several ways such as orthogonalizing the basis, reducing the Rayleigh quotient to Schur form, etc. The similarity transformation effects a change of basis via an invertible matrix $W$,
\begin{equation}\label{eq:similiarity}
AU_kW=U_kW(W^{-1}C_kW)+u_{k+1}c_{k+1}^*W.
\end{equation}
The crux of the Krylov-Schur method for a general matrix $A$ is to obtain a new Rayleigh quotient $W^{-1}C_kW$ with a block of zeros under the diagonal that allows truncating the leading part of the decomposition into a smaller decomposition with the same properties.

Now we are interested in defining a particular subclass of Krylov decomposition and analyzing under which transformations this subclass is closed. For this, we borrow the following notation from~\cite{Shao:2017:PDB},
\begin{equation}\label{eq:notation}
\phi(U,V):=
\begin{bmatrix}
U & \mybar{V} \\
V & \mybar{U} \\
\end{bmatrix},
\end{equation}
where $U$ and $V$ are matrices of the same dimensions. The structure of $\phi$ matrices is preserved under summation, real scaling, complex conjugation, transposition, and matrix multiplication,
\begin{equation}\label{eq:notationmult}
\phi(U_1,V_1)\phi(U_2,V_2)=\phi(U_1U_2+\mybar{V}_1V_2,V_1U_2+\mybar{U}_1V_2).
\end{equation}

\begin{definition}[Structured BSE Krylov decomposition, Type I]\label{def:sbse1}
Let $H=\Hsign\Hhat$ be a definite Bethe--Salpeter matrix, then a structured BSE Krylov decomposition of type I is a relation of the form
\begin{equation}\label{eq:sbse-krylov-1}
H\phi(\ProjBSEX_k,\ProjBSEY_k)
=\phi(\ProjBSEX_k,\ProjBSEY_k)
\begin{bmatrix}
\ProjBSEA_k & \ProjBSEB_k\\
-\ProjBSEB_k & -\ProjBSEA_k
\end{bmatrix}
+
\phi(\ProjBSEx_{k+1},\ProjBSEy_{k+1})
\begin{bmatrix}
a_{k+1}^*&b_{k+1}^*\\
-b_{k+1}^*&-a_{k+1}^*
\end{bmatrix},
\end{equation}
with
\begin{equation}\label{eq:sbse-krylov-1-detail}
A_k=\frac{F_k+K_k}{2},\quad
B_k=\frac{F_k-K_k}{2},\quad
a_{k+1}=\frac{f_{k+1}+k_{k+1}}{2},\quad
b_{k+1}=\frac{f_{k+1}-k_{k+1}}{2},
\end{equation}
for some symmetric positive definite $F_k,K_k\in\mathbb{R}^{k\times k}$ and vectors $f_{k+1},k_{k+1}\in\mathbb{R}^k$,
where the basis of the decomposition (the columns of $\phi([\ProjBSEX_k,\ProjBSEx_{k+1}],[\ProjBSEY_k,\ProjBSEy_{k+1}])$ must be linearly independent. We also require that either $f_{k+1}=0$ or $k_{k+1}=0$, thus $b_{k+1}=\pm a_{k+1}$, which is equivalent to requiring that the last term on the right-hand side (the residual) has rank one.
\end{definition}

\begin{remark}
The Lanczos-type relation~\eqref{eq:projectedbse} satisfies~\cref{def:sbse1} with $F_k=I$, $K_k=T_k$, $a_{k+1}=\beta_ke_k$ and $b_{k+1}=-\beta_ke_k$.
\end{remark}

The decomposition of~\cref{def:sbse1} is a particular instance of a Krylov decomposition. A restarted structure-preserving solver should perform similarity transformations~\eqref{eq:similiarity} with an invertible $W$ that preserves the structure of~\cref{def:sbse1}, as discussed below.

The Rayleigh quotient matrix of the decomposition~\eqref{eq:sbse-krylov-1} is real and has the structure of a linear response eigenvalue problem (LREP) mentioned in \cref{sec:bse}, whose eigenvalues also come in pairs $\{\tilde\lambda_i,-\tilde\lambda_i\}$. The LREP is often solved via an equivalent block matrix whose diagonal blocks are zero and whose off-diagonal blocks are symmetric. This is done via similarity, and is the basis of our definition of the type II decompositions.

\begin{definition}[Structured BSE Krylov decomposition, Type II]\label{def:sbse2}
Let $H=\Hsign\Hhat$ be a definite Bethe--Salpeter matrix, then a structured BSE Krylov decomposition of type II is a relation of the form
\begin{equation}\label{eq:sbse-krylov-2}
H\begin{bmatrix}
U_k & V_k\\
\mybar{U}_k & -\mybar{V}_k
\end{bmatrix}
=
\begin{bmatrix}
U_k & V_k\\
\mybar{U}_k & -\mybar{V}_k
\end{bmatrix}
\begin{bmatrix}
0 & K_k\\F_k&0
\end{bmatrix}
+\begin{bmatrix}
u_{k+1}&v_{k+1}\\\bar{u}_{k+1} & -\bar{v}_{k+1}
\end{bmatrix}
\begin{bmatrix}
0 & k^{*}_{k+1}\\f^{*}_{k+1}&0
\end{bmatrix}
\end{equation}
that has been obtained by applying a similarity transformation with the symmetric orthogonal matrix
$\GruningP=\frac{1}{\sqrt{2}}
\left[\begin{smallmatrix}
I_k&I_k\\I_k&-I_k
\end{smallmatrix}\right]$
to the type I decomposition~\eqref{eq:sbse-krylov-1}, where $F_k$, $K_k$, $f_{k+1}$, $k_{k+1}$ are the same as in~\cref{def:sbse1}.
In particular, $F_k$ and $K_k$ are symmetric positive definite, and either $f_{k+1}=0$ or $k_{k+1}=0$. The two decompositions are equivalent in the sense that they span the same Krylov subspace.
\end{definition}

\begin{proposition}
A structured BSE Krylov decomposition of type II (\cref{def:sbse2}) keeps its structure under block diagonal similarity transformations of the form $\left[\begin{smallmatrix}Q&0\\0&Q^{-*}\end{smallmatrix}\right]$, with $Q\in\mathbb{R}^{k\times k}$ invertible.
\end{proposition}
\begin{proof}
The transformed decomposition would be
\[
H\begin{bmatrix}
\tilde{U}_k & \tilde{V}_k\\
\mybar{\tilde{U}}_k & -\mybar{\tilde{V}}_k
\end{bmatrix}
=
\begin{bmatrix}
\tilde{U}_k & \tilde{V}_k\\
\mybar{\tilde{U}}_k & -\mybar{\tilde{V}}_k
\end{bmatrix}
\begin{bmatrix}
0 & \tilde{K}_k\\\tilde{F}_k&0
\end{bmatrix}
+\begin{bmatrix}
u_{k+1}&v_{k+1}\\\bar{u}_{k+1} & -\bar{v}_{k+1}
\end{bmatrix}
\begin{bmatrix}
0 & \tilde{k}^{*}_{k+1}\\\tilde{f}^{*}_{k+1}&0
\end{bmatrix}
\]
with $\tilde{U}_k:=U_kQ$,
$\tilde{V}_k:=V_kQ^{-*}$,
$\tilde{K}_k:=Q^{-1}K_kQ^{-*}$,
$\tilde{F}_k:=Q^{*}F_kQ$.
$\tilde{k}_{k+1}:=Q^{-1}k_{k+1}$, and
$\tilde{f}_{k+1}:=Q^{*}f_{k+1}$.
Note that both $\tilde{K}_k$ and $\tilde{F}_k$ are symmetric positive definite, since they are obtained from congruence transformations of symmetric positive definite matrices.
\end{proof}

\begin{proposition}\label{prop:sbse1-transform}
A structured BSE Krylov decomposition of type I (\cref{def:sbse1}) keeps its structure under transformations of the form
$\frac{1}{2}\phi(Q+Q^{-*},Q-Q^{-*})$.
\end{proposition}
\begin{proof}
The transformation matrix can be written as
\[
\frac{1}{2}\phi(Q+Q^{-*},Q-Q^{-*})=
\frac{1}{\sqrt{2}}\begin{bmatrix}I_k&I_k\\I_k&-I_k\end{bmatrix}
\begin{bmatrix}Q&0\\0&Q^{-*}\end{bmatrix}
\frac{1}{\sqrt{2}}\begin{bmatrix}I_k&I_k\\I_k&-I_k\end{bmatrix},
\]
which can be seen as a sequence of three transformations, the first of which converts the decomposition from type I to type II, according to \Cref{def:sbse2}, while the second keeps the structure of type II and the third takes the decomposition back to type I.
\end{proof}

\begin{remark}
The decomposition of Shao~\eqref{eq:shao-lanczos} fits the type II decomposition of~\cref{def:sbse2}. 
The corresponding restarted Lanczos method uses orthogonal transformations of block diagonal structure $\left[\begin{smallmatrix}Q&0\\0&Q\end{smallmatrix}\right]$, which are a particular case of the transformation 
$\left[\begin{smallmatrix}Q&0\\0&Q^{-*}\end{smallmatrix}\right]$.

The transformation used in~\eqref{eq:restart-projbse} has also the form
$\left[\begin{smallmatrix}Q&0\\0&Q\end{smallmatrix}\right]$ with orthogonal $Q$, which
is a particular case of the type I structure preserving transformation $\frac{1}{2}\phi(Q+Q^{-*},Q-Q^{-*})$.
\end{remark}

\begin{definition}[Structured BSE Krylov decomposition, Type III]\label{def:sbse3}
Let $H=\Hsign\Hhat$ be a definite Bethe--Salpeter matrix, then a structured BSE Krylov decomposition of type III is a relation of the form
\begin{equation}\label{eq:sbse-krylov-3}
H\begin{bmatrix}
U_k & V_k\\
\mybar{U}_k & -\mybar{V}_k
\end{bmatrix}
=
\begin{bmatrix}
U_k & V_k\\
\mybar{U}_k & -\mybar{V}_k
\end{bmatrix}
\begin{bmatrix}
0 & \hat{K}_k\\\hat{K}_k^*&0
\end{bmatrix}
+\begin{bmatrix}
u_{k+1}&v_{k+1}\\\bar{u}_{k+1} & -\bar{v}_{k+1}
\end{bmatrix}
\begin{bmatrix}
0 & \hat{k}^{*}_{k+1}\\\hat{f}^{*}_{k+1}&0
\end{bmatrix}
\end{equation}
where $\hat{K}_k\in\mathbb{R}^{k\times k}$ is not necessarily symmetric, and either $\hat{f}_{k+1}=0$ or $\hat{k}_{k+1}=0$.
\end{definition}

We note that a decomposition of type III can be obtained by applying a similarity transformation to a decomposition of type II. In particular, 
since $F_k$ is symmetric positive definite, it can be factored as $F_k=\check{F}_k\check{F}_k^*$ and the following similarity transformation can be applied:
\[
\begin{bmatrix}
\check{F}_k^{*} & \\&\check{F}_k^{-1}
\end{bmatrix}
\begin{bmatrix}
0 & K_k\\F_k&0
\end{bmatrix}
\begin{bmatrix}
\check{F}_k^{-*} & \\&\check{F}_k
\end{bmatrix}
=
\begin{bmatrix}
0 & \check{F}^*K_k\check{F}\\I_k&0
\end{bmatrix}.
\]
$K_k$ is also symmetric positive definite, thus we can also factor $\check{F}_k^*K_k\check{F}_k=\hat{K}_k\hat{K}_k^*$, and apply the similarity transformation
\[
\begin{bmatrix}
I_k & \\&\hat{K}_k^*
\end{bmatrix}
\begin{bmatrix}
0 & \check{F}^*K_k\check{F}\\I_k&0
\end{bmatrix}
\begin{bmatrix}
I_k & \\&\hat{K}_k^{-*}
\end{bmatrix}
=
\begin{bmatrix}
0 & \hat{K}_k\\\hat{K}_k^*&0
\end{bmatrix}.
\]
Therefore, the matrix
$\left[\begin{smallmatrix}
\check{F}_k^{-*} & \\&\check{F}_k\hat{K}_k^{-*}
\end{smallmatrix}\right]$
transforms a decomposition of type II into a decomposition of type III.

\begin{proposition}
A structured BSE Krylov decomposition of type III (\cref{def:sbse3}) keeps its structure under transformations of the form
$\left[\begin{smallmatrix}
Q_1&\\&Q_2
\end{smallmatrix}\right]$,
with orthogonal $Q_1, Q_2\in\mathbb{R}^{k\times k}$. 
\end{proposition}

\begin{proof}
The transformed decomposition would be
\[
H\begin{bmatrix}
\tilde{U}_k & \tilde{V}_k\\
\mybar{\tilde{U}}_k & -\mybar{\tilde{V}}_k
\end{bmatrix}
=
\begin{bmatrix}
\tilde{U}_k & \tilde{V}_k\\
\mybar{\tilde{U}}_k & -\mybar{\tilde{V}}_k
\end{bmatrix}
\begin{bmatrix}
0 & \tilde{K}_k\\\tilde{K}^{*}_k&0
\end{bmatrix}
+\begin{bmatrix}
u_{k+1}&v_{k+1}\\\bar{u}_{k+1} & -\bar{v}_{k+1}
\end{bmatrix}
\begin{bmatrix}
0 & \tilde{k}^{*}_{k+1}\\\tilde{f}^{*}_{k+1}&0
\end{bmatrix}
\]
with $\tilde{U}_k:=U_kQ_1$,
$\tilde{V}_k:=V_kQ_2$,
$\tilde{K}_k:=Q_1^{*}K_kQ_2$,
$\tilde{k}_{k+1}:=Q_2^{*}k_{k+1}$, and
$\tilde{f}_{k+1}:=Q_1^{*}f_{k+1}$,
with either $\tilde{k}_{k+1}=0$ or $\tilde{f}_{k+1}=0$.
\end{proof}

By using structured Krylov decompositions of type I, II or III we have a projected matrix which, like the original $H$ matrix, has paired $\{\tilde\lambda_i,-\tilde\lambda_i\}$ eigenvalues and structured eigenvectors. In particular, in the decomposition of type I, the projected matrix is real BSE (LREP), so its eigenvectors retain the structure of \Cref{thm:shao} (although with real eigenvectors). In the case of the decompositions of type II and III, the projected matrix corresponds to an equivalent formulation of the LREP, for which we have a spectral decomposition of the form
\[
\begin{bmatrix}\hat{X}_2&-\hat{X}_2\\\hat{X}_1&\hat{X}_1\end{bmatrix}^*
\begin{bmatrix}0&K\\F&0\end{bmatrix}
\begin{bmatrix}\hat{X}_1&\hat{X}_1\\\hat{X}_2&-\hat{X}_2\end{bmatrix}
=
\begin{bmatrix}\Lambda_+&\\&-\Lambda_+\end{bmatrix}
\]
for the type II, and
\[
\begin{bmatrix}\hat{X}_1&\hat{X}_2\\\hat{X}_1&-\hat{X}_1\end{bmatrix}^*
\begin{bmatrix}0&K\\K^*&0\end{bmatrix}
\begin{bmatrix}\hat{X}_1&\hat{X}_1\\\hat{X}_2&-\hat{X}_2\end{bmatrix}
=
\begin{bmatrix}\Lambda_+&\\&-\Lambda_+\end{bmatrix}
\]
for the type III. The latter relation corresponds to the eigendecomposition of a Jordan-Wielandt matrix~\cite[Th.~1.2.7]{Bjorck:2024:NML}.

\section{Error bounds}\label{sec:error}

We now discuss error bounds for the solutions computed with the proposed algorithms. All methods implicitly rely on a Lanczos recurrence using the $\Hhat$-inner product, $(u,v)_{\Hhat}=v^*\Hhat u$. Therefore, the error analysis should be done by interpreting the original non-Hermitian eigenproblem $Hx=\lambda x$ as a generalized Hermitian-definite eigenvalue problem $\Hsign x=\lambda^{-1}\Hhat x$. It is well known that the conditioning of the inner product matrix $\Hhat$ will affect the error in the computed solution. From this perspective, we could apply standard error bounds for the generalized Hermitian-definite eigenvalue problem, e.g., as described in~\cite[\S 5.7]{Bai:2000:TSA}. Shao and Yang~\cite{Shao:2017:PDB} have adapted those results to the specific structure of definite BSE problems. In particular, they prove that Rayleigh--Ritz-based algorithms that produce small residual norms are backward stable, and this implies forward stability when $\kappa_2(H)$ is not too large.

Suppose $k<n$ Ritz pairs $(\tilde\lambda_i,\tilde x_i)$ have been computed, with $\tilde\lambda_i>0$, $i=1,\dots,k$, and $\tilde\lambda_1\leq\dots\leq\tilde\lambda_k$. Then the following theorem can be proved.

\begin{theorem}[Shao and Yang~\cite{Shao:2017:PDB}]\label{thm:shaoyang}
Let $H=\Hsign\Hhat$ be a definite Bethe--Salpeter matrix and let $\tilde\Lambda_+=\operatorname{diag}\{\tilde\lambda_1,\dots,\tilde\lambda_k\}$ and $\tilde{X}_1,\tilde{X}_2\in\mathbb{C}^{n\times k}$ such that $\tilde x_i=\left[\begin{smallmatrix}\tilde{X}_1\\\tilde{X}_2\end{smallmatrix}\right]e_i$ is the $i$th approximate eigenvector, satisfying
\begin{equation}\label{eq:ritz}
\begin{bmatrix}
\tilde{X}_1 & \mybar{\tilde{X}}_2 \\
\tilde{X}_2 & \mybar{\tilde{X}}_1 \\
\end{bmatrix}^*
\Hhat
\begin{bmatrix}
\tilde{X}_1 & \mybar{\tilde{X}}_2 \\
\tilde{X}_2 & \mybar{\tilde{X}}_1 \\
\end{bmatrix}
=I_{2k},\qquad
\begin{bmatrix}
\tilde{X}_1 & \mybar{\tilde{X}}_2 \\
\tilde{X}_2 & \mybar{\tilde{X}}_1 \\
\end{bmatrix}^*
\Hsign
\begin{bmatrix}
\tilde{X}_1 & \mybar{\tilde{X}}_2 \\
\tilde{X}_2 & \mybar{\tilde{X}}_1 \\
\end{bmatrix}
=
\begin{bmatrix}
\tilde{\Lambda}_+ & 0 \\
0 & -\tilde{\Lambda}_+ \\
\end{bmatrix}^{-1}.
\end{equation}
Define the residual
\begin{equation}\label{eq:residual}
R=H
\begin{bmatrix}
\tilde{X}_1 & \mybar{\tilde{X}}_2 \\
\tilde{X}_2 & \mybar{\tilde{X}}_1 \\
\end{bmatrix}
-
\begin{bmatrix}
\tilde{X}_1 & \mybar{\tilde{X}}_2 \\
\tilde{X}_2 & \mybar{\tilde{X}}_1 \\
\end{bmatrix}
\begin{bmatrix}
\tilde{\Lambda}_+ & 0 \\
0 & -\tilde{\Lambda}_+ \\
\end{bmatrix}.
\end{equation}
Then there exist $k$ positive eigenvalues of $H$, $\lambda_{j_1}\leq\dots\leq\lambda_{j_k}$, such that
\begin{equation}\label{eq:shaoyang}
|\tilde\lambda_i-\lambda_{j_i}|\leq\|H\|_2^{1/2}\|R\|_2,\qquad 1\leq i\leq k.
\end{equation}
\end{theorem}
A tighter bound can also be found in~\cite{Shao:2017:PDB} that takes into account eigenvalue separation.

\section{Details of the implementation in SLEPc}\label{sec:implem}

In this section we discuss a few details related to the implementation of the methods within the SLEPc framework.

\paragraph{User interface}

To solve a Bethe--Salpeter eigenproblem with SLEPc, the application code has to first create the matrices $R$ and $C$, as any valid PETSc \texttt{Mat} type, and then call \texttt{MatCreateBSE(R,C,\&H)}. This helper function will create matrix $H$ of~\eqref{eq:bse1} as a nested matrix (\texttt{MATNEST}) that stores the upper blocks $R$ and $C$ explicitly but handles the bottom blocks in an implicit way.

Once the $2n\times 2n$ matrix $H$ has been created, it is passed to the SLEPc eigensolver as any other matrix. If the problem is solved as a standard non-Hermitian eigenproblem, then the structure will not be taken into account, but if the problem type is set to \texttt{EPS\_BSE} and the default eigensolver is selected then the structure-preserving solvers presented above will be used. All other features such as stopping criterion, convergence monitor, residual check, etc.\ work as usual.

\paragraph{Management of basis vectors}

Every SLEPc eigensolver internally holds a basis of vectors (\texttt{BV}) of the same length as the matrix, i.e., $2n$ in this case. However, the algorithms presented in this paper operate with two bases of length $n$ (e.g., $U_k$ and $V_k$ in the case of \cref{alg:shao-reorth}). For the internal implementation of these solvers, we have added a new operation \texttt{BVGetSplitRows()} that provides a clean way to view the basis of $2n$-vectors as two stacked bases of $n$-vectors, including the case when these vectors are distributed across several MPI processes in a parallel computation.

Once the eigensolver has reached convergence, the basis vectors contain all the necessary information to build the right and left eigenvectors as in \cref{thm:shao}. One remaining thing is how to normalize eigenvectors. We have opted for normalizing them such that $\|x_i\|_2=1$ (which implies that also $\|y_i\|_2=1$), $i=1,\dots,n_\mathrm{ev}$. This is consistent with the normalization in standard non-Hermitian eigenproblems in SLEPc. An alternative would be to normalize them with respect to the norm induced by $\Hhat$ or the pseudo-norm induced by the signature $S$ (since eigenvectors are orthogonal with respect to both inner products), but the application we are interested in (Yambo) expected unit-norm eigenvectors.

\paragraph{Orthogonalization}

As in any Lanczos-type eigensolver, one can propose sophisticated schemes such as partial reorthogonalization in order to avoid loss of orthogonality among the generated Lanczos vectors when Ritz values begin to stabilize. In SLEPc we tend to avoid this type of strategies, due to the increased complexity of implementation and because for moderate restart size full reorthogonalization is usually not much more expensive. In most SLEPc solvers we use Classical Gram-Schmidt with conditional reorthogonalization~\cite{Hernandez:2007:PAE}. However, the methods of \cref{sec:thick} cannot reuse this procedure directly because the formulas are slightly different. In this case, we have opted for a hybrid scheme consisting of (i) local (3-term) orthogonalization, except in the case of restart, as in lines \ref{alg:shao-reorth-reset:line3}--\ref{alg:shao-reorth-reset:line8} of \cref{alg:shao-reorth-reset}, followed by (ii) unconditional full reorthogonalization involving all previously computed Lanczos vectors, as in lines \ref{alg:shao-reorth-reset:line9}--\ref{alg:shao-reorth-reset:line11} of \cref{alg:shao-reorth-reset}. In the experiments of \cref{sec:results}, this scheme has been enough to compute bases that are orthogonal to working precision. However, in difficult problems, depending on the spectrum distribution, it might be necessary to perform an additional full orthogonalization\footnote{Our current implementation does not incorporate a third orthogonalization.}.

\paragraph{Shift-and-invert}

As mentioned in \cref{sec:problemdef}, the discussed methods already approximate the smallest magnitude eigenvalues, which are the relevant ones for applications, without having to apply any kind of spectral transformation. Still, we discuss this possibility here, as it may be useful in case of tightly clustered eigenvalues around zero.

Most SLEPc solvers can be combined with a shift-and-invert technique, which, in the case of standard eigenproblems, will operate with matrix $(A-\sigma I)^{-1}$ to compute eigenvalues closest to $\sigma$. The first thing to note is that in the case of Bethe--Salpeter eigensolvers, we must have $\sigma=0$, otherwise the Bethe--Salpeter structure is lost. It can be easily shown that $H^{-1}$ also has a (definite) Bethe--Salpeter structure. Another comment is that solving linear systems with $H$ using an iterative method such as GMRES is also not viable since again the Bethe--Salpeter structure would not be preserved. We should either build $H^{-1}$ explicitly (possible in the case of dense matrices, but not implemented in SLEPc), or handle $H^{-1}$ implicitly via a Schur complement scheme. The latter can be used in SLEPc via the PETSc preconditioner \texttt{PCFIELDSPLIT} selecting the appropriate options, which would involve two factorizations of size $n$, one for the $R$ block and another one for the Schur complement.

\section{Computational results}\label{sec:results}

We now provide some results about accuracy, convergence and computational performance of the proposed methods. The solvers are available in SLEPc since version 3.22. We use version 3.23, together with PETSc 3.23, built with GNU compilers (version 11.4.0) and MPICH version 4.0. All computational experiments are carried out in complex arithmetic. The computer used has two AMD EPYC GENOA 9354 processors, each of them with 32 physical cores (64 threads) at 3.25 GHz, with 384 GB of main memory (DDR5 at 4800 MHz). The server also includes an AMD Instinct MI210 GPU with 64 GB of HBM2e memory and 6656 stream processors.

\begin{table}
\caption{Description of the test problems used in the computational experiments: type of problem (sparse or dense), numerical precision used, dimension $n$ of the $R$ and $C$ blocks, number of requested eigenvalues \texttt{nev}, first computed (positive) eigenvalue, and average separation of the first \texttt{nev}/2 positive eigenvalues.}
\label{tab:bseproblems}
\begin{minipage}{\columnwidth}
\begin{center}
\begin{revisedblock}
\begin{tabular}{lcccccc}
\hline
name & type & precision & $n$ & \texttt{nev} & first value & sep\\
\hline
\textsf{pentadiag small} & sparse & double & 5000  & 100 & $2.1503397672$ & $3.3\cdot 10^{-5}$\\ 
\textsf{pentadiag large} & sparse & double & 50000 & 100 & $2.1503391439$ & $3.2\cdot 10^{-7}$\\ 
\textsf{CrI$_3$ small}      & dense  & single & 1152  & 100 & $0.0518548376$ & $1.6\cdot 10^{-4}$\\ 
\textsf{CrI$_3$ medium}      & dense  & single & 11520 & 100 & $0.0518094226$ & $1.6\cdot 10^{-4}$\\ 
\textsf{CrI$_3$ large}  & dense  & single & 36288  & 100 & $0.051301926$ & $3.3\cdot 10^{-5}$\\
\textsf{CrI$_3$ verylarge}       & dense  & single & 103680 & 100 & $0.0512986928$ & $3.3\cdot 10^{-5}$\\
\hline
\end{tabular}
\end{revisedblock}
\end{center}
\end{minipage}
\end{table}

\Cref{tab:bseproblems} lists the problems used in the experiments, summarizing some properties and parameters. The first two problems use a synthetic matrix created to exercise the solvers with sparse matrices, but it does not correspond to any real application. The \revisedtext{third and fourth}{last four} problems come from a simulation with the \yambo code, where the matrices are dense. In the case of \yambo, we use single precision arithmetic, which is sufficient for the application requirements and provides substantial memory savings in the case of large problems. Here is a short description of the problems:
\begin{itemize}
\item In the \textsf{pentadiag} test, both $R$ and $C$ blocks are Toeplitz matrices with only a few nonzero diagonals. More precisely, $R = \mathrm{pentadiag}(a,b,c,\bar{b},\bar{a})$ and $C = \mathrm{tridiag}(b,d,b)$. In our tests we set the values $a=-0.1+0.2\text{i}$, $b=1+0.5\text{i}$, $c=4.5$, and $d=2+0.2\text{i}$. The tolerance in this case is set to $10^{-8}$, the default for double precision.
\item \textsf{CrI$_3$} (chromium triiodide) is a dense matrix from a \yambo simulation. The size of the blocks $R$ and $C$ is $n_k n_c n_v$ where $n_k$ is related to the size of the grid, $n_c$ is the number of unoccupied states and $n_v$ is the number of occupied states. We work with \revisedtext{two}{four} cases of different size\revisedtext{, with the same value of $n_k$, but different}{. For the first two matrices, $n_k=36$, and for the last two, $n_k=324$. They also differ in the} band intervals, a parameter which modifies $n_c$ and $n_v$. The band interval is 31-40 for the \textsf{small} case\revisedtext{}{, 31-48 for the \textsf{large} case} and 23-48 for the \revisedtext{\textsf{large} case}{\textsf{medium} and \textsf{verylarge} cases}. The tolerance for this problem is set to $10^{-5}$, the default for single precision.
\end{itemize}

In all cases, 100 eigenvalues are computed (\texttt{nev}).  The number of column vectors (\texttt{ncv}, which corresponds to $k$ in the algorithms of \cref{sec:thick}) is set to twice the number of requested eigenpairs that will be computed by the iterative solver. The methods of \cref{sec:thick} effectively compute half of the requested eigenpairs, which in this case are 50 eigenvalues from the positive side of the spectrum and their eigenvectors, while their negative counterparts are obtained a posteriori. Therefore an \texttt{ncv} of 100 is used. For the non-Hermitian method, used for comparison, \texttt{ncv} is set to 200, because all 100 requested eigenvalues are computed by the iterative method.
\begin{revisedblock}
All the methods performed locking of converged eigenpairs (see \cref{sec:eigentriplets}).
\end{revisedblock}

\begin{table}
\caption{Results of the computational experiments in CPU with one MPI process, and in GPU. The solvers are compared in terms of the number of restarts (its), execution time, largest relative residual and bi-orthogonality of the right and left eigenvectors. The non-Hermitian method computes only the right eigenvectors, the left eigenvectors are built as in \cref{thm:shao} to check bi-orthogonality.}
\label{tab:bsecpuresults}
\begin{minipage}{\columnwidth}
\begin{center}
\begin{revisedblock}
\begin{tabular}{lccccccc}
\hline
name & exec & solver & its & time & residual & bi-orthog\\
\hline
\textsf{pentadiag small} & CPU & \shao          & 152  & 6.48      & $2.60\cdot10^{-9}$ & $1.34\cdot10^{-14}$ \\
\textsf{pentadiag small} & CPU & \gruning       & 151  & 9.9       & $5.00\cdot10^{-9}$ & $2.73\cdot10^{-10}$ \\
\textsf{pentadiag small} & CPU & \projectedbse  & 152  & 6.88      & $2.60\cdot10^{-9}$ & $1.69\cdot10^{-14}$ \\
\textsf{pentadiag small} & CPU & \nhe & 151  & 38.52     & $9.86\cdot10^{-9}$ & $5.33\cdot10^{-6}$  \\ 
\textsf{pentadiag large} & CPU & \shao          & 6520 & 4,778     & $2.60\cdot10^{-9}$ & $8.78\cdot10^{-14}$ \\
\textsf{pentadiag large} & CPU & \gruning       & 6613 & 5,600     & $5.72\cdot10^{-9}$ & $6.23\cdot10^{-10}$ \\
\textsf{pentadiag large} & CPU & \projectedbse  & 6537 & 5,014     & $2.59\cdot10^{-9}$ & $1.05\cdot10^{-13}$ \\
\textsf{pentadiag large} & CPU & \nhe & 4064 & 15,418    & $9.16\cdot10^{-9}$ & $2.99\cdot10^{-2}$  \\ 
\textsf{CrI$_3$ small}   & CPU & \shao          & 16   & 0.60      & $7.74\cdot10^{-5}$ & $1.97\cdot10^{-6}$  \\
\textsf{CrI$_3$ small}   & CPU & \gruning       & 17   & 0.75      & $5.09\cdot10^{-6}$ & $3.42\cdot10^{-6}$  \\
\textsf{CrI$_3$ small}   & CPU & \projectedbse  & 16   & 0.65      & $9.86\cdot10^{-5}$ & $3.63\cdot10^{-6}$  \\
\textsf{CrI$_3$ small}   & CPU & \nhe & 15   & 2.00      & $2.69\cdot10^{-5}$ & $9.20\cdot10^{-1}$  \\ 
\textsf{CrI$_3$ medium}   & CPU & \shao          & 48   & 394       & $7.94\cdot10^{-5}$ & $3.76\cdot10^{-6}$  \\
\textsf{CrI$_3$ medium}   & CPU & \gruning       & 50   & 416       & $1.93\cdot10^{-5}$ & $4.93\cdot10^{-6}$  \\
\textsf{CrI$_3$ medium}   & CPU & \projectedbse  & 51   & 403       & $7.88\cdot10^{-5}$ & $7.61\cdot10^{-6}$  \\
\textsf{CrI$_3$ medium}   & CPU & \nhe & 47   & 876       & $6.54\cdot10^{-5}$ & $7.79\cdot10^{-1}$  \\ 
\textsf{CrI$_3$ large}   & CPU & \shao          &  56  &  4,785      & $7.49\cdot10^{-5}$ & $2.07\cdot10^{-5}$  \\
\textsf{CrI$_3$ large}   & CPU & \gruning       &  68  &  6,003      & $6.17\cdot10^{-5}$ & $1.41\cdot10^{-5}$  \\
\textsf{CrI$_3$ large}   & CPU & \projectedbse  &  53  &  4,690      & $7.28\cdot10^{-5}$ & $1.36\cdot10^{-5}$  \\
\textsf{CrI$_3$ large}   & CPU & \nhe           &  73  &  13,480     & $7.54\cdot10^{-5}$ & $9.17\cdot10^{-1}$  \\
\textsf{CrI$_3$ verylarge}   & CPU & \shao          & 88   & 100,890    & $2.12\cdot10^{-4}$ & $5.45\cdot10^{-5}$  \\
\textsf{CrI$_3$ verylarge}   & CPU & \gruning       & 98   & 113,910    & $3.60\cdot10^{-4}$ & $5.12\cdot10^{-5}$  \\
\textsf{CrI$_3$ verylarge}   & CPU & \projectedbse  & 73   &  85,970    & $1.83\cdot10^{-4}$ & $5.34\cdot10^{-5}$  \\
\textsf{CrI$_3$ verylarge}   & CPU & \nhe           & 90   & 234,400    & $9.73\cdot10^{-5}$ & $7.57\cdot10^{-1}$  \\
\hline
\textsf{pentadiag small} & GPU & \shao          & 152  & 5.84   & $2.60\cdot10^{-9}$ & $1.55\cdot10^{-14}$ \\
\textsf{pentadiag small} & GPU & \gruning       & 151  & 6.96   & $5.00\cdot10^{-9}$ & $2.73\cdot10^{-10}$ \\
\textsf{pentadiag small} & GPU & \projectedbse  & 152  & 9.00  & $2.60\cdot10^{-9}$ & $1.50\cdot10^{-14}$ \\
\textsf{pentadiag small} & GPU & \nhe & 150  & 15.99  & $9.71\cdot10^{-9}$ & $5.60\cdot10^{-6}$  \\ 
\textsf{pentadiag large} & GPU & \shao          & 6628 & 579    & $2.65\cdot10^{-9}$ & $1.11\cdot10^{-13}$ \\
\textsf{pentadiag large} & GPU & \gruning       & 6483 & 625    & $5.56\cdot10^{-9}$ & $6.16\cdot10^{-10}$ \\
\textsf{pentadiag large} & GPU & \projectedbse  & 6572 & 717    & $2.61\cdot10^{-9}$ & $1.15\cdot10^{-13}$ \\
\textsf{pentadiag large} & GPU & \nhe & 3280 & 542    & $9.57\cdot10^{-9}$ & $4.84\cdot10^{-2}$  \\ 
\textsf{CrI$_3$ small}   & GPU & \shao          & 16   & 0.69   & $9.70\cdot10^{-5}$ & $2.03\cdot10^{-6}$  \\
\textsf{CrI$_3$ small}   & GPU & \gruning       & 17   & 0.86   & $5.65\cdot10^{-6}$ & $3.13\cdot10^{-6}$  \\
\textsf{CrI$_3$ small}   & GPU & \projectedbse  & 15   & 0.92   & $1.01\cdot10^{-4}$ & $4.30\cdot10^{-6}$  \\
\textsf{CrI$_3$ small}   & GPU & \nhe & 21   & 2.14   & $2.73\cdot10^{-5}$ & $9.53\cdot10^{-1}$  \\ 
\textsf{CrI$_3$ medium}   & GPU & \shao          & 50   & 9.23   & $8.39\cdot10^{-5}$ & $6.51\cdot10^{-6}$  \\
\textsf{CrI$_3$ medium}   & GPU & \gruning       & 50   & 10.33  & $1.15\cdot10^{-5}$ & $6.33\cdot10^{-6}$  \\
\textsf{CrI$_3$ medium}   & GPU & \projectedbse  & 51   & 10.73  & $8.20\cdot10^{-5}$ & $6.53\cdot10^{-6}$  \\
\textsf{CrI$_3$ medium}   & GPU & \nhe & 92   & 30.94  & $8.84\cdot10^{-5}$ & $5.27\cdot10^{-1}$  \\ 
\textsf{CrI$_3$ large}   & GPU & \shao          &  56  & 84   & $7.94\cdot10^{-5}$ & $3.34\cdot10^{-6}$  \\
\textsf{CrI$_3$ large}   & GPU & \gruning       &  67  & 105  & $9.75\cdot10^{-6}$ & $4.91\cdot10^{-6}$  \\
\textsf{CrI$_3$ large}   & GPU & \projectedbse  &  61  & 91   & $6.65\cdot10^{-5}$ & $8.12\cdot10^{-6}$  \\
\textsf{CrI$_3$ large}   & GPU & \nhe           &  69  & 219  & $7.40\cdot10^{-5}$ & $6.11\cdot10^{-1}$  \\
\hline
\end{tabular}
\end{revisedblock}
\end{center}
\end{minipage}
\end{table}

\Cref{tab:bsecpuresults} presents results of our three new solvers (\shao, \gruning, \projectedbse) as well as non-Hermitian Krylov-Schur for comparison, with executions in both CPU and GPU. The number of restarts required to compute 100 eigenvalues clearly indicates that the \textsf{pentadiag} problem is much more difficult compared to \textsf{CrI$_3$}. The justification is that eigenvalues of \textsf{pentadiag} are closer to each other, especially for large dimension. The average separation of the wanted eigenvalues is shown in the last column of \cref{tab:bseproblems}. \revisedtext{All three structure-preserving methods take approximately the same number of restarts, while t}{In some cases, t}he non-Hermitian solver \revisedtext{sometimes }{}needs less iterations, maybe because the size of the subspace (\texttt{ncv}) is larger. Still, the execution time of the non-Hermitian solver is in general significantly larger than that of the structured solvers. Among the structured solvers, \gruning is slightly slower than the other two in CPU, probably due to having to read a second vector basis from memory.\revisedtext{ When running on GPU, it is \projectedbse that needs more time, which indicates some inefficiencies probably related to synchronizing data between the GPU and CPU memory.}{}

In terms of accuracy, all solutions computed by all solvers have a relative residual smaller than the tolerance in double precision calculations, while in single precision the residual is slightly larger than the tolerance in some cases. \Cref{tab:bsecpuresults} shows the maximum relative residual of all computed eigentriplets $(\tilde\lambda_i,\tilde{x}_i,\tilde{y}_i)$, $i=1,\dots,\texttt{nev}$, evaluated as
\begin{equation}
\max\{\|H\tilde{x}_i-\tilde\lambda_i\tilde{x}_i\|_2,\|\tilde{y}_i^*H-\tilde{y}_i^*\tilde\lambda_i\|_2\}/|\tilde\lambda_i|.
\end{equation}
Remember that eigenvectors are normalized such that $\|\tilde{x}_i\|_2=1$ and $\|\tilde{y}_i\|_2=1$. The table also shows a quantification of the bi-orthogonality of computed eigenvectors, evaluated as $\|\tilde{Y}^*\tilde{X}-\Delta\|_\mathrm{max}$, where $\Delta=\texttt{diag}(\tilde{y}_i^*\tilde{x}_i)$. We can see that the structured eigensolvers enforce bi-orthogonality up to \revisedtext{}{or near} the prescribed tolerance, while eigenvectors computed by the non-Hermitian solver fail to satisfy this property, as expected.
Another comment is that the structured eigensolvers return real eigenvalues, while in the case of the non-Hermitian solver the eigenvalues will have an imaginary part that can be as large as $1.4\cdot 10^{-10}$ in \textsf{pentadiag} and \revisedtext{$1.4\cdot 10^{-7}$}{$1.6\cdot 10^{-7}$} in \textsf{CrI$_3$}.

\begin{table}%
\caption{Maximum residual, over the sampled restarts, of the Lanczos decomposition equality.}
\label{tab:resid-lanczos}
\begin{minipage}{\columnwidth}
\begin{center}
\begin{revisedblock}
\begin{tabular}{lcccc}
\hline
Problem                   & \shao    & \gruning & \projectedbse & \nhe\\
\hline
\textsf{pentadiag small}  & $6.05\cdot10^{-8}$ & $1.87\cdot10^{-8}$ & $4.28\cdot10^{-8}$ & $3.70\cdot10^{-8}$ \\
\textsf{pentadiag large}  & $1.05\cdot10^{-7}$ & $2.92\cdot10^{-8}$ & $8.45\cdot10^{-8}$ & $4.02\cdot10^{-8}$ \\
\textsf{CrI$_3$ small}    & $6.04\cdot10^{-6}$ & $2.39\cdot10^{-6}$ & $3.67\cdot10^{-6}$ & $1.66\cdot10^{-6}$ \\
\textsf{CrI$_3$ medium}    & $5.66\cdot10^{-6}$ & $9.73\cdot10^{-6}$ & $8.32\cdot10^{-6}$ & $3.46\cdot10^{-6}$ \\
\end{tabular}
\end{revisedblock}
\end{center}
\end{minipage}
\end{table}

\begin{table}%
\caption{Maximum departure from orthogonality of the basis, over the sampled restarts.}
\label{tab:resid-orthog}
\begin{minipage}{\columnwidth}
\begin{center}
\begin{revisedblock}
\begin{tabular}{lcccc}
\hline
Problem                   & \shao    & \gruning & \projectedbse & \nhe\\
\hline
\textsf{pentadiag small}  & $1.02\cdot10^{-13}$ & $1.74\cdot10^{-9}$ & $4.44\cdot10^{-14}$ & $3.05\cdot10^{-13}$\\
\textsf{pentadiag large}  & $6.52\cdot10^{-13}$ & $2.95\cdot10^{-9}$ & $3.45\cdot10^{-13}$ & $6.52\cdot10^{-12}$\\
\textsf{CrI$_3$ small}    & $1.97\cdot10^{-5}$ & $5.84\cdot10^{-6}$ & $1.36\cdot10^{-5}$ & $3.60\cdot10^{-5}$\\
\textsf{CrI$_3$ medium}    & $4.05\cdot10^{-5}$ & $1.60\cdot10^{-5}$ & $2.12\cdot10^{-5}$ & $8.97\cdot10^{-5}$\\
\end{tabular}
\end{revisedblock}
\end{center}
\end{minipage}
\end{table}

\Cref{tab:resid-lanczos,tab:resid-orthog} provide information about the accuracy of the computed Lanczos process in each method \revisedtext{and for each problem}{for the first 4 problems}.
In order to obtain this information, we have \revisedtext{used a Matlab implementation of the methods, in which}{computed} the residuals of the Lanczos decompositions and the level of basis orthogonality \revisedtext{are obtained}{} every few restarts. In particular, we compute the residuals of the truncated
Lanczos decompositions \eqref{eq:shao-lanczos-trunc}, \eqref{eq:gruning-pre-truncated} or \eqref{eq:projectedbse-restart}, depending of the method, e.g., in the first method, the residual is given by the matrix norm
\[
\left\|H
\begin{bmatrix}
\hat{U}_r            & \hat{V}_r\\
\mybar{\hat{U}}_r & -\mybar{\hat{V}}_r
\end{bmatrix}
-
\begin{bmatrix}
\hat{U}_r            & \hat{V}_r\\
\mybar{\hat{U}}_r & -\mybar{\hat{V}}_r
\end{bmatrix}
\begin{bmatrix}
0  &D_r\\
I_r & 0
\end{bmatrix}
\revisedtext{}{-}
\begin{bmatrix}
u_{k+1}\\\bar{u}_{k+1}
\end{bmatrix}
\begin{bmatrix}
0  &b^*
\end{bmatrix}\right\|,
\]
where $r$ is the number of columns kept in the restart. Similarly, orthogonality is measured taking into account the corresponding basis orthogonality equalities \eqref{eq:shao-ortog}, \eqref{eq:gruning-orthog} or \eqref{eq:projectedbse-ortog}, e.g., in the first case,
\[
\left\|
\begin{bmatrix}
V_r        & U_r\\
\mybar{V}_r&-\mybar{U}_r
\end{bmatrix}^*
\begin{bmatrix}
U_r         & V_r\\
\mybar{U}_r & -\mybar{V}_r
\end{bmatrix}
-
2I_{2r}
\right\|.
\]
\begin{revisedblock}
In the case of the \nhe method, the Lanczos residual and orthogonality deviation can be computed respectively as
\[
\begin{gathered}
	\|H \mathring{U}_r - \mathring{U}_r \mathring{T}_r - \mathring{u}_{k+1} \mathring{b}^*\|,\\
	\|\mathring{U}_r^*\mathring{U}_r-I_r\|,
\end{gathered}
\]
where the basis vectors  $\mathring{U}_r$ have no structure, and $\mathring{T}_r$ is a complex upper triangular matrix resulting from the reduction of the projected matrix to Schur form.
Note that, in order to keep the same number of vectors, in this case $r$ should be twice as large as in the structured methods.
\end{revisedblock}

\revisedtext{We observed that the residuals increased slowly with the number of restarts, the maximum values being obtained at the last iterations.}{}
\Cref{tab:resid-lanczos} presents the maximum residual of the Lanczos decomposition for each method, while \Cref{tab:resid-orthog} shows the maximum departure from orthogonality.
We can see that the different methods are comparable in both respects\revisedtext{, although \gruning presents higher residuals for the \textsf{pentadiag large} matrix}.
\revisedtext{}{Due to space constraints, we do not show graphs with the evolution of the residuals and the orthogonality deviation as the methods converge. However, we observed that both the residuals and the orthogonality deviation increased slowly with the number of restarts, the maximum values being obtained at the last iterations. This is something to be expected, affecting also the \nhe method, based on the standard inner product. However, since the proposed methods use the $\hat{H}$ inner product, the issue may get worse in case of an ill-conditioned $\hat{H}$.
We also observed that locking converged eigenpairs makes the Lanczos residual grow up to near the order of the prescribed tolerance, which is also to be expected. In the case of \gruning, locking has a similar effect on the orthogonality deviation, while it does not affect orthogonality of other methods. That accounts for the larger values of orthogonality deviation in \cref{tab:resid-orthog} for \gruning on the \textsf{pentadiag} problems.}

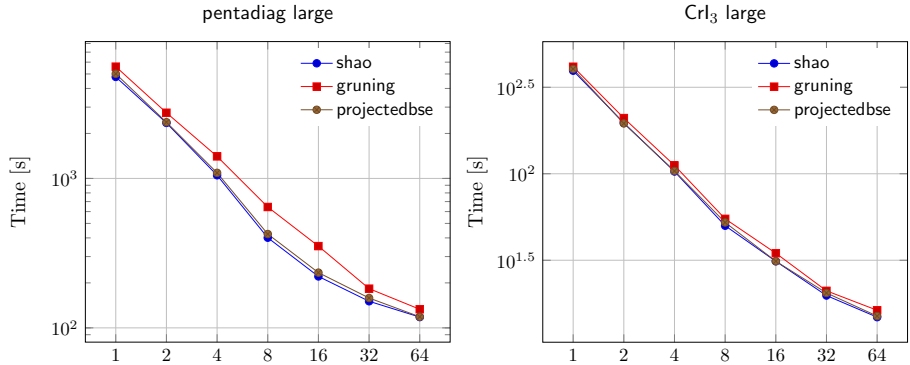
\begin{figure}
\label{fig:parallel}
\centering
\begin{tikzpicture}[scale=0.7] 
  \begin{loglogaxis}[
  title={\textsf{pentadiag large}},
  ylabel={Time [s]},
  grid=major,
  log basis x=2,
  xtick={1,2,4,8,16,32,64},
  xticklabels={1,2,4,8,16,32,64},
  ticklabel style={font=\small},
  legend style={draw=none, legend columns=1,font=\small,cells={anchor=west}}
  ]
  \addplot coordinates { 
  	(  1, 4777.70)
  	(  2, 2351.30)
  	(  4, 1049.40)
  	(  8, 401.09)
  	( 16, 220.91)
  	( 32, 150.66)
  	( 64, 118.07)
  };
  \addplot coordinates { 
  	(  1, 5599.70)
  	(  2, 2752.20)
  	(  4, 1403.80)
  	(  8, 643.97)
  	( 16, 352.50)
  	( 32, 183.07)
  	( 64, 133.21)
  };
  \addplot coordinates { 
  	(  1, 5013.90)
  	(  2, 2378.30)
  	(  4, 1090.00)
  	(  8, 424.52)
  	( 16, 234.22)
  	( 32, 158.43)
  	( 64, 118.66)
  };
	\legend{\shao, \gruning, \projectedbse}
  \end{loglogaxis}
\end{tikzpicture}
\begin{tikzpicture}[scale=0.7] 
  \begin{loglogaxis}[
  title={\textsf{CrI$_3$ large}},
  ylabel={Time [s]},
  grid=major,
  log basis x=2,
  xtick={1,2,4,8,16,32,64},
  xticklabels={1,2,4,8,16,32,64},
  ticklabel style={font=\small},
  legend style={draw=none, legend columns=1,font=\small,cells={anchor=west}}
  ]
  \addplot coordinates { 
  	(  1, 4784.8)
  	(  2, 2399.2)
  	(  4, 1193.8)
  	(  8, 650.10)
  	( 16, 378.79)
  	( 32, 196.97)
  	( 64, 158.52)
  };
  \addplot coordinates { 
  	(  1, 6003.0)
  	(  2, 3212.6)
  	(  4, 1589.1)
  	(  8, 760.34)
  	( 16, 470.65)
  	( 32, 245.74)
  	( 64, 211.98)
  };
  \addplot coordinates { 
  	(  1, 4689.6)
  	(  2, 2440.7)
  	(  4, 1284.8)
  	(  8, 649.77)
  	( 16, 405.54)
  	( 32, 220.16)
  	( 64, 159.15)
  };
	\legend{\shao, \gruning, \projectedbse}
  \end{loglogaxis}
\end{tikzpicture}
\begin{tikzpicture}[scale=0.7] 
  \begin{loglogaxis}[
  title={\textsf{pentadiag large}},
  ylabel={Speedup},
  grid=major,
  log basis x=2,
  xtick={1,2,4,8,16,32,64},
  xticklabels={1,2,4,8,16,32,64},
  ticklabel style={font=\small},
  legend pos = south east,
  legend style={draw=none, legend columns=1,font=\small,cells={anchor=west}}
  ]
  \addplot coordinates { 
  	(  1, 4777.70/4777.70)
  	(  2, 4777.70/2351.30)
  	(  4, 4777.70/1049.40)
  	(  8, 4777.70/401.09)
  	( 16, 4777.70/220.91)
  	( 32, 4777.70/150.66)
  	( 64, 4777.70/118.07)
  };
  \addplot coordinates { 
  	(  1, 5599.70/5599.70)
  	(  2, 5599.70/2752.20)
  	(  4, 5599.70/1403.80)
  	(  8, 5599.70/643.97)
  	( 16, 5599.70/352.50)
  	( 32, 5599.70/183.07)
  	( 64, 5599.70/133.21)
  };
  \addplot coordinates { 
  	(  1, 5013.90/5013.90)
  	(  2, 5013.90/2378.30)
  	(  4, 5013.90/1090.00)
  	(  8, 5013.90/424.52)
  	( 16, 5013.90/234.22)
  	( 32, 5013.90/158.43)
  	( 64, 5013.90/118.66)
  };
	\legend{\shao, \gruning, \projectedbse}
  \end{loglogaxis}
\end{tikzpicture}
\begin{tikzpicture}[scale=0.7] 
  \begin{loglogaxis}[
  title={\textsf{CrI$_3$ large}},
  ylabel={Speedup},
  grid=major,
  log basis x=2,
  xtick={1,2,4,8,16,32,64},
  xticklabels={1,2,4,8,16,32,64},
  ticklabel style={font=\small},
  legend pos = south east,
  legend style={draw=none, legend columns=1,font=\small,cells={anchor=west}}
  ]
  \addplot coordinates { 
  	(  1, 4784.8/4784.8)
  	(  2, 4784.8/2399.2)
  	(  4, 4784.8/1193.8)
  	(  8, 4784.8/650.10)
  	( 16, 4784.8/378.79)
  	( 32, 4784.8/196.97)
  	( 64, 4784.8/158.52)
  };
  \addplot coordinates { 
  	(  1, 6003.0/6003.0)
  	(  2, 6003.0/3212.6)
  	(  4, 6003.0/1589.1)
  	(  8, 6003.0/760.34)
  	( 16, 6003.0/470.65)
  	( 32, 6003.0/245.74)
  	( 64, 6003.0/211.98)
  };
  \addplot coordinates { 
  	(  1, 4689.6/4689.6)
  	(  2, 4689.6/2440.7)
  	(  4, 4689.6/1284.8)
  	(  8, 4689.6/649.77)
  	( 16, 4689.6/405.54)
  	( 32, 4689.6/220.16)
  	( 64, 4689.6/159.15)
  };
	\legend{\shao, \gruning, \projectedbse}
  \end{loglogaxis}
\end{tikzpicture}
\caption{
Parallel execution time and speedup on CPU for the \textsf{pentadiag large} test on the left, and \textsf{CrI$_3$ large} test on the right, for an increasing number of MPI processes.
}
\end{figure}

Finally, we evaluate the parallel scalability of the new solvers. \Cref{fig:parallel} plots execution times on CPU of the new eigensolvers for \revisedtext{the two large}{two of the largest} test cases when increasing the number of MPI processes. As expected, scalability is a bit worse in the case of dense matrices. A more comprehensive performance study \revisedtext{with \yambo matrices up to size 103680}{} can be found in~\cite{Milev:2026:SHP}.

\section{Concluding remarks}\label{sec:concl}

We have developed three new Lanczos-type eigensolvers for the Bethe--Salpeter eigenvalue problem. The methods rely on previously proposed structure-preserving Lanczos recurrences, but we have incorporated all the details required for a robust and efficient practical implementation. In particular, all the details related to thick restart and reorthogonalization have been worked out. The result is production-quality solvers, incorporated in the SLEPc library, that can be readily used by applications such as \yambo. According to~\cite{Milev:2026:SHP}, the performance of the eigenvalue solution in \yambo simulations has improved in an order of magnitude, compared to the previous solution approach (two-sided non-Hermitian Krylov-Schur with the same value of the \texttt{ncv} parameter). This will enable \yambo users to perform more detailed analyses involving larger matrices with a more effective use of supercomputer allocated time. The results in \cref{sec:results} also show that the computed solution is more accurate as far as the level of bi-orthogonality of left and right eigenvectors is concerned.

In this paper we have assumed that the number of wanted eigenvalues is in the order of hundreds or less. However, some of the application articles using the \yambo code that were cited in~\cref{sec:intro} compute a few thousand eigenvalues. In that case, the performance of our solvers will likely degrade due to the increased cost of full orthogonalization. An alternative would be to use some kind of semi-orthogonalization, but the robustness of this approach might be compromised when combined with thick restart. As an ongoing effort, we are considering spectrum slicing strategies, where eigenvalues are computed in chunks. In particular, we are working out the details of a polynomial filtering method adapted to the Bethe--Salpeter matrix.

\begin{acknowledgements}
We would like to thank Davide Sangalli for providing us with the \yambo test matrices used in \cref{sec:results}, and for fruitful discussion and interaction that led to improvements in the interface between \yambo and SLEPc.
We also wish to thank the anonymous referees for their helpful suggestions and constructive remarks.
\end{acknowledgements}


\end{document}